\documentclass[DIV=14,letterpaper]{scrartcl}

\date{}

\usepackage{tikz}
\usetikzlibrary{calc}
\usetikzlibrary{matrix}
\usepackage{amsmath,amssymb,amsfonts,amsthm,subfig}
\usepackage{hyperref}                 

\usepackage{graphicx}
\allowdisplaybreaks

\DeclareOldFontCommand{\bf}{\normalfont\bfseries}{\mathbf}
\DeclareOldFontCommand{\it}{\normalfont\itshape}{\mathit}

\title{A Third-Order Moving Mesh Cell-Centered Scheme for One-Dimensional Elastic-Plastic Flows}

\author{Jun-Bo Cheng%
\thanks{Laboratory of Computational Physics, Institute of Applied Physics and
Computational Mathematics, Beijing, China ({\em cheng\_junbo@iapcm.ac.cn})}
\and
Weizhang~Huang%
\thanks{Department of Mathematics, the University of Kansas, Lawrence, KS 66045, U.S.A.
({\em whuang@ku.edu})}
\and
Song Jiang%
\thanks{Laboratory of Computational Physics, Institute of Applied Physics and
Computational Mathematics, Beijing, China ({\em jiang@iapcm.ac.cn})}
\and
Baolin Tian%
\thanks{Laboratory of Computational Physics, Institute of Applied Physics and
Computational Mathematics, Beijing, China ({\em tian\_baolin@iapcm.ac.cn})}
}

\begin{document}
\vskip 1cm
\maketitle

\begin{abstract}
A third-order moving mesh cell-centered scheme without the remapping of physical variables
is developed for the numerical solution of one-dimensional elastic-plastic flows with the Mie-Gr\"{u}neisen equation of state,  the Wilkins constitutive model, and the von Mises yielding criterion. The scheme combines the Lagrangian method with the MMPDE moving mesh method and adaptively moves the mesh to better resolve shock and other types of waves while preventing the mesh from crossing and tangling. It can be viewed as a direct arbitrarily Lagrangian-Eulerian method but can also be degenerated to a purely Lagrangian scheme. It treats the relative velocity of the fluid with respect to the mesh as constant in time between time steps, which allows high-order approximation of
free boundaries. A time dependent scaling is used in the monitor function to avoid possible sudden movement of the mesh points due to the creation or diminishing of shock and rarefaction waves or the steepening of those waves. A two-rarefaction Riemann solver with elastic waves is employed to compute the Godunov values of the density, pressure, velocity, and deviatoric stress at cell interfaces. Numerical results are presented for three examples. The third-order convergence of the scheme and its ability to concentrate mesh points around shock and elastic rarefaction waves are demonstrated. The obtained numerical results are in good agreement with those in literature. The new scheme is also shown to be more accurate in resolving shock and rarefaction waves than an existing third-order cell-centered Lagrangian scheme.
\end{abstract}



\noindent{\bf Abbreviated title.}
A third-order moving mesh scheme for elastic-plastic flows

\section{Introduction}

We consider the numerical simulation of one-dimensional (1D) elastic-plastic solid
problems with the isotropic elastic-plastic model initially developed by Wilkins
\cite{wilkins} where a perfectly elastic material is characterized by Hooke's law in the sense that
an incremental strain results in an incremental stress and the elastic limit is described by
the von Mises yielding criterion.

Numerical methods for elastic-plastic flows can be
roughly classified into three categories. The first is staggered Lagrangian schemes,
first developed by Wilkins \cite{wilkins}, which discretize
the equations of momentum and specific internal energy on a staggered mesh
and use artificial viscosity to suppress spurious numerical oscillations around shock waves.
The second is Eulerian methods, e.g., see \cite{euler4,euler1,euler5,euler2,euler3},
which are often used for hyperelastic models that can be formulated in the conservative form.
The third, which is considered in this work,
is cell-centered Lagrangian schemes that have recently attracted considerable attention
from researchers (e.g., see \cite{cc1,HLLCE,trrse,cc2,cc3,cc4}).
They are conservative and can resolve shock waves
without using artificial viscosity and be used for both hyperelastic and hypoelastic models.
However, the major disadvantage of these schemes is that the mesh, which moves at the flow velocity,
can easily cross over or become tangling.

Arbitrary Lagrangian-Eulerian (ALE) methods, first introduced by Hirt, Amsden, and Cook \cite{Hirt},
have been designed to overcome this disadvantage. Usually, ALE methods contain three phases,
the pure Lagrangian phase, a rezone phase, and a remapping phase. 
Different from general ALE methods, direct ALE methods have been implemented directly without the need for interpolation (i.e. remapping) of physical variables between different meshes.
For example, Luo et al. \cite{ale-luo} introduce second-order
direct ALE HLLCE and Godunov schemes to solve the problems of multi-material flows while
high-order direct ALE finite volume and discontinuous Galerkin (DG) schemes are proposed in
\cite{cale-fluid,cale-hpr} and \cite{aledg1,ale-rdg}, respectively.
Unlike their traditional counterparts, which can readily handle free boundary domains in the Lagrangian
phase and allow the mesh to move arbitrarily in the rezone phase, these direct ALE methods
are not considered to deal with free boundary problems where the boundary moves at the flow velocity
(which is part of the solutions to be sought) and therefore cannot be purely Lagrangian \cite{cale-hpr,aledg1}.
Direct ALE methods have been successfully applied to hyperelastic-platic flows in \cite{cale-hpr}.

Interestingly, ALE methods form a special type of moving mesh methods.
In the past, a number of other moving mesh methods have been developed;  e.g., see
\cite{BHJ05a, BW06, mmpde2, HR99, LTZ01, MM81}. The interested reader is also
referred to the books/review articles \cite{Bai94a,Baines-2011,BHR09,HR11,Tan05} and references therein.
In this work, we consider combining an cell-centered Lagrangian scheme with
an error-based moving mesh method which is known to be more robust in preventing
the mesh from crossing and tangling.
We will use the so-called MMPDE (moving mesh partial differential equation) moving mesh
method \cite{mmpde1, mmpde2, HR11} which is capable of concentrating
mesh points according to a user-supplied monitor function and preventing the mesh from becoming singular
in one and multiple dimensions \cite{HK2015}. It also has the advantages of being relatively simple,
compact, and easy to implement and working for both convex and concave domains \cite{mmpde4}.
By combining the Lagrangian and MMPDE methods, we hope that the new method will inherit
the advantages of both methods. Moreover, although we focus on
1D problems in the current work, we hope that the new method can be used for multi-dimensional problems.

More specifically, we will develop a third-order moving
mesh cell-centered scheme (MMCC) for 1D elastic-plastic flows.
The scheme can be viewed as a new direct ALE method.
Indeed, like the existing direct ALE methods \cite{cale-fluid,cale-hpr,ale-luo,aledg1,ale-rdg}, MMCC scheme
does not need the remapping of physical variables between different meshes. But it also possesses
a unique feature that other direct ALE methods do not have, that is,
it can be degenerated to a purely Lagrangian scheme. This is done by introducing
the relative velocity of the fluid with respect to the moving mesh as a new variable $w$.
When $w = 0$, the MMCC scheme reduces to a purely Lagrangian scheme (CCL) proposed
by Cheng et al. \cite{trrse}.
Moreover, the introduction of $w$ makes it relatively easy to simulate free boundary problems for which we can simply
set $w=0$ on the boundary (so the boundary nodes move at the Lagrangian velocity).
Furthermore, the mesh velocity is commonly treated in existing moving mesh methods
as constant in time between time steps, which typically results in low-order approximation
of free boundaries (e.g, see \cite{Huang2014b}). MMCC treats $w$ as constant in time instead,
and this new treatment allows high-order approximation of free boundaries.
In fact, like CCL, MMCC is also third-order accurate (and conservative and essentially non-oscillatory).
In addition, it has the ability to concentrate mesh points around shock waves and elastic rarefaction waves.

An outline of the paper is as follows. The governing equations of elastic-plastic flows
are described in Section~\ref{SEC:pde}.
The MMCC scheme is presented in Section~\ref{SEC:MMCC},
followed by numerical examples in Section~\ref{SEC:numerics}.
The conclusions are drawn in Section~\ref{SEC:conclusions}.

\section{Governing equations}
\label{SEC:pde}

The 1D governing equations of the isotropic elastic-plastic model of \cite{wilkins}
for elastic-plastic flows read as
\begin{equation}   \label{e1}
      \frac{\partial}{\partial t} {\bf U} + \frac{\partial}{\partial x} {\bf F}({\bf U})
      = 0, \quad  x \in \Omega,\quad t>0
\end{equation}
where $\Omega$ is an interval,
\begin{equation}   \label{e2}
U{\text{ = }}\left[ {\begin{array}{*{20}{c}}
  \rho  \\
  {\rho u} \\
  {\rho E}
\end{array}} \right],
\qquad F{\text{ = }}\left[ {\begin{array}{*{20}{c}}
  {\rho u} \\
  {\rho {u^2} - {\sigma _x}} \\
  {\left(\rho E - {\sigma _x}\right)u}
\end{array}} \right],
\end{equation}
$\rho$ is the density, $u$ is the velocity in the $x$-direction,
$E=e+\frac{1}{2}u^{2}$ is the total energy density, $e$ is
the specific internal energy, and $\sigma_x$ is the Cauchy stress
which is related to the hydrostatic pressure ($p$) and the deviatoric stress ($s_{xx}$) as
\begin{equation}
\label{e3}
     \sigma _x = -p + s_{xx} .
\end{equation}
Following the convention in mechanics, we use $\sigma_x$ and $s_{xx}$ here to denote
the components of the corresponding variables in $x$ direction; this is different from the mathematical convention
where subscription $x$ typically denotes the differentiation with respect to $x$.
The pressure is related to the internal energy and the density through the equation of state (EOS).
In this work we use the Mie-Gr\"{u}neisen EOS,
\begin{equation}   \label{e4a}
     p(\rho, e)=\rho _0 a_0^2 \frac{(\eta-1)\left[\eta -
\Gamma_{0}(\eta-1)/2\right]}{\left[\eta-s(\eta-1)\right]^{2}}
 + \rho _0 \Gamma _0 e,
\end{equation}
where $\eta = {\rho}/{\rho_{0}}$ and $\rho_0$, \textbf{$a_0$}, $s$, and
$\Gamma_0$ are constant parameters.
The elastic energy is not included in the above model as it is negligible in many
practical engineering problems (e.g., see \cite{cc3}).

The relation between the deviatoric stress and the strain is described by
Hooke's law. It reads as
\begin{equation}   \label{e6a}
\dot{s}_{xx} = 2\mu \Big( \dot{\varepsilon}_x -
\frac{1}{3}\frac{\dot{V}}{V}\Big),
\end{equation}
where $\mu$ is the shear modulus, $V$ is volume, the dot denotes the material time derivative,
and
\[
\dot{\varepsilon}_x=\frac{\partial u}{\partial x},
\quad
\frac{\dot{V}}{V}=\frac{\partial u}{\partial x}.
\]
Combining the above equations we get
\begin{equation}   \label{e10}
\frac{\partial s_{xx}}{\partial t}+u\frac{\partial
s_{xx}}{\partial x} = \frac{4 \mu}{3} \frac{\partial u}{\partial x}.
\end{equation}

The von Mises yielding condition is used here to describe the
elastic limit. In one dimension, it takes the form
\begin{equation}   \label{e11}
|s_{xx}| \leq \frac{2}{3}Y_{0},
\end{equation}
where $Y_{0}$ is the yield strength of the material in simple tension.
This condition is enforced in the numerical solution of (\ref{e1}) and (\ref{e10});  cf. (\ref{e64}).

For the numerical examples we consider in Section~\ref{SEC:numerics}, free/fixed boundaries with
periodic, Dirichlet, and/or wall boundary conditions are used, depending on the specific setting
of each example. Thus, it is important that any scheme designed for (\ref{e1}) and (\ref{e10}) can
approximate free boundaries with high-order accuracy.

\section{A third-order moving mesh cell-centered scheme}
\label{SEC:MMCC}

\subsection{The overall procedure}
\label{SEC:overall}

With the MMCC method, the governing equations (\ref{e1}) and (\ref{e10}) are discretized
on a moving mesh
\begin{equation}
x_{\frac 1 2}(t) < x_{\frac 3 2}(t) < \cdots < x_{N-\frac 1 2}(t) < x_{N+\frac 1 2}(t),
\label{mesh-1}
\end{equation}
where $N$ is the number of cells in the mesh.
The mesh is moved using a combination of the Lagrangian and
MMPDE moving mesh velocities \cite{mmpde1, mmpde2, HR11}. This is different from most
of the existing moving mesh methods (e.g., see \cite{BHR09,HR11})
where one of the velocities, but not both, is used.
By combining the velocities, we hope that the new method can inherit the advantages
of both, especially being able to concentrate mesh points around shock waves and some
other places of interest (for instance, the region of large gradient of $s_{xx}$ in our current situation)
while preventing the mesh from crossing and tangling.

Another main difference between MMCC and the existing moving mesh methods
lies in the numerical treatment of the mesh speed ($\dot{x}$).
It is common practice (e.g., see \cite{HR11}) that $\dot{x}$ is approximated to be constant in time between
time steps. A disadvantage of this is that the trajectories of the boundary points are linear in time,
which gives a second-order approximation to free boundaries
(that move at the flow velocity in the current situation).
To allow higher-order approximations, we propose to treat the relative
flow velocity with respect to the mesh, $w = u - \dot{x}$, as constant between time steps instead.
With this, we can assure that $w$ be zero and thus $\dot{x} = u$ on the free boundary.
Thus, the scheme is purely Lagrangian on the free boundary. Moreover,
a high-order approximation to the location of the free boundary can be obtained by
integrating $\dot{x} = u$ (see the detail in Section~\ref{SEC:CC}).

The steps of the MMCC method are given in the following and will be elaborated in the subsequent subsections.

\begin{itemize}
\item[(i)] Given the cell average $\bar{Q}$ of $Q = (\rho, \rho u, \rho E, s_{xx})^{T}$, the Lagrangian
	velocity $u_{i+\frac{1}{2}}^{n}$, $i=0,\ldots,N$ is evaluated
	using the two-rarefaction Riemann solver with elastic wave (TRRSE) \cite{trrse}.
Notice that (\ref{e10}) is not in conservative form. TRRSE is specially designed for the nonconservative
system (\ref{e1}) and (\ref{e10}) in primitive variables by assuming that all nonlinear waves are
continuous rarefaction waves and solving the equations for the generalized Riemann invariants
across rarefactions; see \cite{trrse} for the detail.
Like a two-rarefaction Riemann solver (TRRS) in \cite{toro} for pure fluids, TRRSE has proven
robust and accurate for elastic-plastic flows \cite{trrse}.

\item[(ii)] The Lagrangian coordinate of the mesh points are computed as
\begin{equation}\label{e109}
  \widehat{x}_{i+\frac{1}{2},L}^{n+1}=x_{i+\frac{1}{2}}^{n} + u_{i+\frac{1}{2}}^{n}\Delta t^{n-{\frac{1}{2}}},
  \quad  i=0,\ldots,N
\end{equation}
where $\Delta t^{n-\frac{1}{2}}=t^{n}-t^{n-1}$.

\item[(iii)] (Section~\ref{SEC:mmpde}) Based on $\widehat{x}_{i+\frac{1}{2},L}^{n+1}$, $i = 0, \ldots, N$,
	the MMPDE method is used to generate a new mesh $\widehat{x}_{i+\frac{1}{2}}^{n+1}$, $i=0,\ldots,N$
	with boundary conditions $\widehat{x}_{\frac{1}{2}}^{n+1} = \widehat{x}_{\frac{1}{2},L}^{n+1}$
	and $\widehat{x}_{N+\frac{1}{2}}^{n+1} = \widehat{x}_{N+\frac{1}{2},L}^{n+1}$.

\item[(iv)] The relative velocity is computed as
 \begin{equation}\label{e110}
  w_{i+\frac{1}{2}} = u_{i+\frac{1}{2}}^{n} - \frac{ \widehat{x}_{i+\frac{1}{2}}^{n+1}-{x}_{i+\frac{1}{2}}^{n}}
  {\Delta t^{n-\frac{1}{2}}},\quad i=0,\ldots,N
\end{equation}
and is treated as constant in time from $t_n$ to $t_{n+1}$.
Notice that $w_{\frac{1}{2}}=w_{N+\frac{1}{2}}=0$ and thus the boundary mesh nodes
move at the flow velocity.

\item[(v)] (Section~\ref{SEC:CFL}) The new time step size $\Delta t^{n}$ is estimated
using the CFL condition. Set $t^{n+1}=t^{n}+\Delta t^{n}$.

\item[(vi)] (Section~\ref{SEC:CC}) The third-order Runge-Kutta method and third-order finite volume numerical
	fluxes are used to compute $\overline{Q}_{i}^{n+1}$, $i=1,\ldots,N$, and $x_{i+\frac{1}{2}}^{n+1}$, $i=0,\ldots,N$.
	The MMCC method is expected to be of third order in both time and space.
\end{itemize}

\subsection{Mesh movement by the MMPDE method}
\label{SEC:mmpde}

The MMPDE method is a variational method that generates the mesh through
a coordinate transformation $x = x(\xi, t): (0,1) \to (a, b)$ (with
$a = \widehat{x}_{\frac{1}{2},L}^{n+1}$ and $b=\widehat{x}_{N+\frac{1}{2},L}^{n+1}$
being fixed for this time step), viz.,
\[
x_{i+\frac 1 2}(t) = x(\hat{\xi}_{i+\frac 1 2}, t), \quad i = 0, ..., N
\]
where $\hat{\xi}_{i+\frac 1 2} = i/N$, $i = 0, \ldots, N$ is a uniform mesh on $(0,1)$.
The inverse of the coordinate transformation, $\xi = \xi(x,t)$,  is governed by an MMPDE defined as
the gradient flow equation of a meshing functional. We use the functional
\begin{equation}
I[\xi] = \frac 1 2 \int_a^b \frac{1}{M} \left (\frac{\partial \xi}{\partial x}\right )^2 d x,
\label{fun-1}
\end{equation}
where $M$ is the monitor function that is defined to measure difficulty in approximating the physical solution.
The Euler-Lagrange equation of the functional is
\[
- \frac{\partial }{\partial x} \left (\frac{1}{M} \frac{\partial \xi}{\partial x}\right ) = 0,
\]
which is mathematically equivalent to the well-known equidistribution principle,
\[
M \frac{\partial x}{\partial \xi} = \sigma,  \qquad \sigma = \int_a^b M d x .
\]
Then the MMPDE is defined as
\begin{equation}
\frac{\partial \xi}{\partial t} = \frac{M}{\tau} \frac{\partial }{\partial x} \left (\frac{1}{M} \frac{\partial \xi}{\partial x}\right ),
\label{mmpde-1}
\end{equation}
where $\tau > 0$ is a parameter used to adjust the response time scale of the mesh movement
to the changes in $M$. The smaller $\tau$ is, the more quickly the mesh movement responds to the changes in $M$.
We take $\tau = 0.01$ in our computation.

Recall that our goal is to generate the new mesh $\widehat{x}_{i+\frac{1}{2}}^{n+1}$, $i=0,\ldots,N$
from the Lagrangian mesh $\widehat{x}_{i+\frac{1}{2},L}^{n+1}$, $i = 0, \ldots, N$.
To this end, we discretize (\ref{mmpde-1}) on the Lagrangian mesh using central finite differences in space
and the backward Euler in time. We have
\begin{align}
& \frac{\xi_{i+\frac{1}{2}}^{n+1}-\hat{\xi}_{i+\frac{1}{2}}}{\Delta t}
= \frac{M_{i+\frac{1}{2}}}{\tau (\widehat{x}_{i+1,L}^{n+1}-\widehat{x}_{i,L}^{n+1})}\left[
\frac{1}{M_{i+1}} \frac{(\xi_{i+\frac{3}{2}}^{n+1}-\xi_{i+\frac{1}{2}}^{n+1})}{(\widehat{x}_{i+\frac{3}{2},L}^{n+1}-\widehat{x}_{i+\frac{1}{2},L}^{n+1})}
-\frac{1}{M_{i}} \frac{(\xi_{i+\frac{1}{2}}^{n+1}-\xi_{i-\frac{1}{2}}^{n+1})}{(\widehat{x}_{i+\frac{1}{2},L}^{n+1}-\widehat{x}_{i-\frac{1}{2},L}^{n+1})} \right],
\label{mmpde-2}
\end{align}
where $M_{i}=\frac{1}{2}(M_{i+\frac{1}{2}}+M_{i-\frac{1}{2}})$,
$\widehat{x}_{i,L}=\frac{1}{2}(\widehat{x}_{i+\frac{1}{2},L}+\widehat{x}_{i-\frac{1}{2},L})$, and
the uniform computational mesh has been taken as the initial mesh for the integration.
The above equation,  together with the boundary conditions,
\[
\xi_{\frac{1}{2}}^{n+1} = 0, \quad \xi_{N+\frac{1}{2}}^{n+1} = 1,
\]
forms a linear algebraic system that can be solved for the new computational mesh
$\xi^{n+1}_{i+\frac 1 2}$, $i = 0, \ldots, N$.
This mesh and the Lagrangian mesh define a correspondence relation, denoted by
\[
x = \Phi_h(\xi)\quad \mbox{ or }\quad \widehat{x}^n_{i+\frac 1 2,L} = \Phi_h(\xi^{n+1}_{i+\frac 1 2}),
\quad i = 0, \ldots, N.
\]
Then, the new physical mesh is defined as
\[
\widehat{x}^{n+1}_{i+\frac 1 2} = \Phi_h(\hat{\xi}_{i+\frac 1 2}),\quad i = 0, ..., N
\]
which can be obtained through linear interpolation.

We now consider the definition of the monitor function.
The entropy is commonly used for adaptive mesh simulation of shock waves.
However, in the current situation the complexity of the EOS and the presence of the deviatoric stress
make it difficult to compute the entropy exactly. For this reason, we use the density and deviatoric stress,
which have jumps around shock waves and elastic limits, respectively.
First, the derivatives of $\rho$ and $s_{xx}$ are computed,
\[
\left(\frac{\partial \rho}{\partial x}\right)_{i+\frac{1}{2}}
=\frac{\overline{\rho}_{i+1,L}^{n+1}-\overline{\rho}_{i,L}^{n+1}}{\widehat{x}_{i+1,L}^{n+1}-\widehat{x}_{i,L}^{n+1}},
\qquad \left(\frac{\partial s_{xx}}{\partial x}\right)_{i+\frac{1}{2}}
=\frac{(\overline{s_{xx}})_{i+1}^{n}-(\overline{s_{xx}})_{i}^{n}}{{x}_{i+1}^{n}-{x}_{i}^{n}},
\]
where
\[
\overline{\rho}_{i+1,L}^{n+1}=\frac{\overline{\rho}_{i+1}^{n}\left(x_{i+\frac{1}{2}}^{n}-x_{i-\frac{1}{2}}^{n}\right)}
{\widehat{x}_{i+\frac{1}{2},L}^{n+1}-\widehat{x}_{i-\frac{1}{2},L}^{n+1}} .
\]
Notice that the value of $s_{xx}$ at $t_n$ has been used since it is not constant along a particle trajectory.
Having obtained the derivatives, we define
\begin{equation}
\label{mf2}
\widehat{M}_{i+\frac{1}{2}}^{(n)}=\sqrt{1+\left(\frac{\partial \rho}{\partial x}\right)_{i+\frac{1}{2}}^2+\alpha \left(\frac{\partial s_{xx}}{\partial x}\right)_{i+\frac{1}{2}}^2},
\end{equation}
where
\begin{equation}
\label{magnify}
\alpha = \left[\frac{\max\limits_{i}{|\left(\frac{\partial \rho}{\partial x}\right)_{i+\frac{1}{2}}|}}
{\max\limits_{i}|\left(\frac{\partial s_{xx}}{\partial x}\right)_{i+\frac{1}{2}}|}\right]^{2} .
\end{equation}
The scaling factor $\alpha$ is used in the above since the magnitude of the jumps in the density
is much larger than those in the deviatoric stress in general. This scaling is necessary for concentrating
mesh points around the elastic limit.

Since jumps such as shock waves can disappear or be created,  the ratio of the maximum value
to the minimum value of the monitor function defined in (\ref{mf2}) can change dramatically over time.
To avoid these dramatic changes, we further scale the monitor function using
the current and past ratios of the maximum and minimum values of the monitor function and define
\begin{equation}\label{mf1}
M_{i+\frac{1}{2}}=\widehat{M}_{Min}^{(n)}+\frac{M_{critic}-\widehat{M}_{Min}^{(n)}}
{\widehat{M}_{Max}^{(n)}-\widehat{M}_{Min}^{(n)}}
\left(\widehat{M}_{i+\frac{1}{2}}^{(n)}-\widehat{M}_{Min}^{(n)}\right),
\end{equation}
where
\begin{align*}
& \widehat{M}_{Min}^{(n)}=\min_{i}\left(\widehat{M}_{i+\frac{1}{2}}^{(n)}\right),
\quad \widehat{M}_{Max}^{(n)}=\max_{i}\left(\widehat{M}_{i+\frac{1}{2}}^{(n)}\right),
\\
& M_{critic} = \min\left(10,\max_{k=1,\ldots,n}\frac{\widehat{M}_{Max}^{(k)}}{\widehat{M}_{Min}^{(k)}}\right)
\widehat{M}_{Min}^{(n)} .
\end{align*}

The so-defined monitor function is further smoothed in space in order to produce a mesh that is also smooth
in space. We use
\[
\frac{1}{4} M_{i+\frac{3}{2}} + \frac{1}{2} M_{i+\frac{1}{2}} + \frac{1}{4} M_{i-\frac{1}{2}}
\; \to \; M_{i+\frac{1}{2}}
\]
for the interior points and a similar formula for the boundary points. This scheme is applied
$N/40$ sweeps every time a new monitor function is computed.

\subsection{A third-order cell-centered scheme on moving meshes}
\label{SEC:CC}

We now describe a third-order cell-centered discretization for (\ref{e1}) and (\ref{e10})
under appropriate initial and boundary conditions on the moving mesh (\ref{mesh-1}).
Denote $I_i=(x_{i-\frac{1}{2}},x_{i+\frac{1}{2}})$, $\Delta x_i=x_{i+\frac{1}{2}}-x_{i-\frac{1}{2}}$, and
the center of $I_i$ by $x_{i}$.  A semi-discrete finite volume discretization
of the conservative equation (\ref{e1}) over the cell $I_i$ is given by
\begin{equation}
\label{e92d}
\frac{d(\mathbf{\overline{U} }_{i}\Delta
x_{i})}{dt}=-\left(\mathbf{F}_{i+\frac{1}{2}}-\mathbf{F}_{i-\frac{1}{2}} \right),
\end{equation}
where $\mathbf{\overline{U} }_{i}$ is the numerical approximation to the average
of $U$ over $I_i$ and
$\mathbf{F}_{i+ \frac{1}{2}}$ is the numerical flux at the cell boundary $x_{i+ \frac{1}{2}}$
(and $\mathbf{F}_{i- \frac{1}{2}}$ is the numerical flux at $x_{i- \frac{1}{2}}$).
It is defined as
\begin{equation}   \label{e92e}
 \mathbf{F}_{i+\frac{1}{2}} = F(\mathbf{Q}_{i+\frac{1}{2}})
 =  \begin{bmatrix} \rho_{i+\frac{1}{2}} w_{i+\frac{1}{2}} \\
                      p_{i+\frac{1}{2}}-(s_{xx})_{i+\frac{1}{2}} + \rho_{i+\frac{1}{2}} u_{i+\frac{1}{2}} w_{i+\frac{1}{2}}\\
                      \left(p_{i+\frac{1}{2}}-(s_{xx})_{i+\frac{1}{2}}\right) u_{i+\frac{1}{2}}
                      + \rho_{i+\frac{1}{2}} E_{i+\frac{1}{2}} w_{i+\frac{1}{2}} \end{bmatrix},
\end{equation}
where $\mathbf{Q}_{i+\frac{1}{2}}$ denotes the Godunov value
of $\mathbf{Q} = (\rho, \rho u, \rho E, s_{xx})^T$ to the Riemann problem with the left
and right states, $\mathbf{Q}_{i+\frac{1}{2}}^{-}$ and $\mathbf{Q}_{i+\frac{1}{2}}^{+}$ (which
are the left and right limits of $\mathbf{Q}$ at $x_{i+\frac{1}{2}}$, respectively).
Note that mesh movement has been taken into consideration and $\dot{x}$ has been replaced by
$\dot{x} = u - w$. In the following, we explain how to compute $\mathbf{Q}_{i+\frac{1}{2}}$ and
$\mathbf{Q}_{i+\frac{1}{2}}^{\pm}$.

The limits $\mathbf{Q}_{i+\frac{1}{2}}^{\pm}$ are computed using the
third-order WENO reconstruction (see Jiang and Shu \cite{WENO-JS} for the detail)
in local characteristic variables. We first transform the conservative
variables to the characteristic variables and then use the
third-order WENO reconstruction method to reconstruct each component
of the characteristic variables. The final values are obtained by
transforming back to the conservative variables.
Notice that the Jacobian matrix of (\ref{e1}) and (\ref{e10}) with respect to $\mathbf{Q}$ is
\begin{equation}   \label{e14a}
\mathbf{J}= \begin{bmatrix}
   0 & 1 & 0 & 0  \\
   { - u^2  + \frac{{\partial p}}{{\partial \rho }} + \Gamma \left( {\frac{{u^2 }}{2} - e} \right)}
   & {u\left( {2 - \Gamma } \right)} & \Gamma  & { - 1}  \\
   {\left( {\Gamma \left( {\frac{{u^2 }}{2} - e} \right) - e + \frac{{\sigma _x }}{\rho }
   + \frac{{\partial p}}{{\partial \rho }}} \right)u} & { - \Gamma u^2  - \frac{{\sigma _x }}{\rho } + e}
   & {\left( {1 + \Gamma } \right)u} & { - u}  \\
   {\frac{4}{3}\mu \frac{u}{\rho }} & { - \frac{4}{3}\mu \frac{1}{\rho }} & 0 & u
\end{bmatrix},
\end{equation}
where $ \Gamma=\frac{\Gamma_0\rho_{0}}{\rho}$.
Denote the left and right eigenvectors of $\mathbf{J}$ by $\mathbf{L}$ and $\mathbf{R}$.
Then, the procedure for computing $\mathbf{Q}_{i+\frac{1}{2}}^{\pm}$ is as follows.
\begin{itemize}
\item[(1)] Evaluate
\[
\overline{\mathbf{q}}_{i}=\mathbf{L}(\overline{\mathbf{Q}}_{i}^{n}) \cdot \overline{\mathbf{Q}}_{i}^{n}, \quad i=-2,\ldots,N+2.
\]
\item[(2)] Given $\overline{\mathbf{q}}_{i}$, $i=-2,\ldots,N+2$, use the
third-order WENO reconstruction method to approximate $\mathbf{q}$ at the cell interface, $\mathbf{q}_{i+\frac{1}{2}}^{\pm}$, $i=0,\ldots,N$.
\item[(3)] Evaluate
\[\mathbf{Q}_{i+\frac{1}{2}}^{-}=\mathbf{R}(\overline{\mathbf{Q}}_{i}^{n})  \cdot  {\mathbf{q}}_{i+\frac{1}{2}}^{-},
\qquad \mathbf{Q}_{i-\frac{1}{2}}^{+}=\mathbf{R}(\overline{\mathbf{Q}}_{i}^{n})  \cdot  {\mathbf{q}}_{i-\frac{1}{2}}^{+},
\quad i=0,\ldots,N.
\]
\end{itemize}

\begin{figure}[htp]
\begin{center}
\includegraphics[width=6cm]{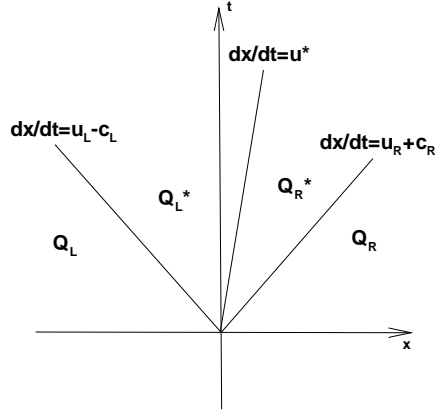}
\end{center}
\caption{A sketch of the solution structure of the Riemann problem (\ref{e108w}).} \label{struc}
\end{figure}

Having obtained $\mathbf{Q}_{i+\frac{1}{2}}^{\pm}$ and letting $\mathbf{Q}_{L} = \mathbf{Q}_{i+\frac{1}{2}}^{-}$
and $\mathbf{Q}_{R} = \mathbf{Q}_{i+\frac{1}{2}}^{+}$, we have the Riemann problem
\begin{equation}   \label{e108w}
 \begin{cases}
\frac{\partial}{\partial t} {\bf U} + \frac{\partial}{\partial x} {\bf F}({\bf U}) = 0,  \\
\frac{\partial}{\partial t} s_{xx} + u \frac{\partial}{\partial x} s_{xx} = \frac{3}{4}\mu\frac{\partial}{\partial x} u,\\
{\bf Q}(x,0) = \left\{
       \begin{array}{lll}
            {\bf Q}_{L}, & \mbox{   for   } & x <  x_{i+\frac{1}{2}} \\[2mm]
            {\bf Q}_{R}, & \mbox{   for   } & x >  x_{i+\frac{1}{2}}.
       \end{array}\right.
\end{cases}
\end{equation}
The structure of the solution to this Riemann problem is sketched in Fig.~\ref{struc}.
It consists of three waves corresponding to three eigenvalues ($\frac{dx}{dt}=u_L-c_L,u^{\ast},u_R+c_R$), where
\[
c = \sqrt {{a^2}  - \frac{{\rho _0 }}{\rho^2 }\Gamma _0 s_{xx}+
\frac{4}{3}\frac{\mu }{\rho }},\quad
a^2  = \frac{{\partial p}}{{\partial \rho }} + \frac{p}{{\rho
^2}}\frac{{\partial p}}{{\partial e}} = a_0^2 \frac{{\partial
f}}{{\partial \eta }} + \frac{p}{{\rho ^2 }}\rho _0 \Gamma _0,
\quad
f=\frac{(\eta-1)\left[\eta - \Gamma_{0}(\eta-1)/2\right]}{\left[\eta-s(\eta-1)\right]^{2}} .
\]
$\frac{dx}{dt}=u^{\ast}$ corresponds to the contact wave while $\frac{dx}{dt}=u_L-c_L$ and
$\frac{dx}{dt}=u_R+c_R$ correspond to the left-going and right-going waves, respectively.
These three waves separate four constant states. From
left to right, they are $\mathbf{Q}_{L}$ (left data state), $\mathbf{Q}_{L}^{\ast}$,
$\mathbf{Q}_{R}^{\ast}$, and $\mathbf{Q}_{R}$ (right data state).
TRRSE \cite{trrse} is used to obtain $\mathbf{Q}_{L}^{\ast}$ and $\mathbf{Q}_{R}^{\ast}$.
Have obtained $\mathbf{Q}_{L}$, $\mathbf{Q}_{L}^{\ast}$,
$\mathbf{Q}_{R}^{\ast}$, and $\mathbf{Q}_{R}$, we define
the Godunov values of $\mathbf{Q}$ as

\begin{equation}                       \label{e19.37}
\mathbf{Q}_{i+\frac{1}{2}}=
\left\{ {\begin{array}{*{20}c}
   {\mathbf{Q}_{_{L} }          }, & {for\;\;  \dot{x} \leq u_{L}-c_{L} }  \\[1mm]
   {\mathbf{Q}_{_{L} }^ {\ast}  }, & {for\;\;  u_{L}-c_{L} < \dot{x} \leq u^ {\ast} }  \\[1mm]
    {\mathbf{Q}_{_{R} }^ {\ast} }, & {for\;\;  u^ {\ast} <  \dot{x} \leq  u_{R}+c_{R}}\\[1mm]
   {\mathbf{Q}_{_{R} }          }, & {for\;\;  \dot{x} >  u_{R}+c_{R}} .
\end{array}} \right.
\end{equation}
Note that the effects of mesh movement (at the speed of $\dot{x}=u-w$) has been taken into consideration.

The Godunov values of $\rho$, $u$, $p$, and $s_{xx}$ can then be obtained from $\mathbf{Q}_{i+\frac{1}{2}}$.
After that, the numerical flux can be computed using (\ref{e92e}) for the conservative equations (\ref{e1}).
It is noted that we have $u^ {\ast} = u_L^ {\ast} = u_R^ {\ast}$ by the construction of TRRSE \cite{trrse}.

We now discuss the spatial discretization of the constitutive equation.
Notice that the equation of the constitute model (\ref{e10}) can be written as
\begin{equation}   \label{e115w}
\left(\frac{{d{s_{xx}}}}{{dt}}\right)_{\dot{x}} = \frac{{4\mu }}{3}\frac{{\partial u}}{{\partial x}}-w\frac{\partial s_{xx}}{\partial x},
\end{equation}
where the left-hand side stands for the time derivative of $s_{xx}$ along the mesh trajectories.
We discretize $\frac{\partial u}{\partial x}$ in cell $I_i$ by
\begin{equation}   \label{e107e.1}
\frac{\partial u}{\partial x} \approx \frac{u_{i+\frac{1}{2}}-u_{i-\frac{1}{2}}}{x_{i+\frac{1}{2}}-x_{i-\frac{1}{2}}}.
\end{equation}
Moreover,
$w\frac{\partial s_{xx}}{\partial x}$ can be discretized as
\begin{align}
 w\frac{\partial s_{xx}}{\partial x} & = \frac{\partial w s_{xx}}{\partial x} -s_{xx} \frac{\partial w }{\partial x}
 \notag \\
 & \approx\frac{w_{i+\frac{1}{2}}s_{xx,i+\frac{1}{2}}
-w_{i-\frac{1}{2}}s_{xx,i-\frac{1}{2}}}{\Delta x_{i}}-(\overline{s_{xx}})_{i}\frac{w_{i+\frac{1}{2}}
-w_{i-\frac{1}{2}}}{\Delta x_{i}} .
\label{e107e.2}
\end{align}
Combining the above approximations, we get
\begin{align}
\frac{ds_{xx}}{dt} & = \frac{{4\mu}}{3}\frac{{u_{i + \frac{1}{2}} -
u_{i - \frac{1}{2}}}}{{x_{i + \frac{1}{2}} - x_{i - \frac{1}{2}}}}-\frac{w_{i+\frac{1}{2}}s_{xx,i+\frac{1}{2}}
-w_{i-\frac{1}{2}}s_{xx,i-\frac{1}{2}}}{\Delta x_{i}}
+(\overline{s_{xx}})_{i}\frac{w_{i+\frac{1}{2}}
-w_{i-\frac{1}{2}}}{\Delta x_{i}} ,
\label{e108}
\end{align}
where $u_{i\pm\frac{1}{2}}$ and $s_{xx,i\pm\frac{1}{2}}$ are the Godunov values of $u$ and
$s_{xx}$ at the left and right boundaries of $I_i$, respectively.

The semi-discrete scheme (\ref{e92d}) and (\ref{e108}) and the von Mises yielding condition (\ref{e11})
are integrated using a third-order TVD Runge-Kutta method. Notice that
the mesh changes with time,  the node location and the cell size need to be
updated at each Runge-Kutta stage. Moreover, the relative velocity $w$ is kept constant
in time for this time step. The scheme reads as follows.

\begin{itemize}
\item[] {\bf Stage 1.}
$$  x_{i+\frac{1}{2}}^{(1)}=x_{i+\frac{1}{2}}^{(0)}+\Delta t^{n} \left(u_{i+\frac{1}{2}}^{n}-w_{i+\frac{1}{2}}\right), $$
$$ \Delta x_{i}^{(1)}=x_{i+\frac{1}{2}}^{(1)}-x_{i-\frac{1}{2}}^{(1)}, $$
$$
\Delta x_{i}^{(1)}\overline{\mathbf{U}}_{i}^{(1)}=\Delta x_{i}^{(0)}\overline{\mathbf{U}}_{i}^{(0)}+\Delta t^{n}
\mathcal{L}\left(\overline{\mathbf{U}}_{i}^{(0)},(\overline{s_{xx}})_{i}^{(0)},x_{i+\frac{1}{2}}^{(0)},w_{i+\frac{1}{2}}\right), $$
$$ (\overline{\widehat{s}_{xx}})_{i}^{(1)}=(\overline{s_{xx}})_{i}^{(0)}+\Delta
t^{n}\Theta\left({u}_{i+\frac{1}{2}}^{(0)},(\overline{s_{xx}})_{i}^{(0)},x_{i+\frac{1}{2}}^{(0)},w_{i+\frac{1}{2}}\right),  $$
$$ (\overline{s_{xx}})_{i}^{(1)}=\Upsilon((\overline{\widehat{s}_{xx}})_{i}^{(1)}); $$
\item[] {\bf Stage 2.}
$$ x_{i+\frac{1}{2}}^{(2)}=\frac{3}{4}x_{i+\frac{1}{2}}^{(0)}+\frac{1}{4}\left[x_{i+\frac{1}{2}}^{(1)}
+\Delta t^{n} \left(u_{i+\frac{1}{2}}^{(1)}-w_{i+\frac{1}{2}}\right)\right], $$
$$ \Delta x_{i}^{(2)}=x_{i+\frac{1}{2}}^{(2)}-x_{i-\frac{1}{2}}^{(2)},  $$
$$ \Delta x_{i}^{(2)}\overline{\mathbf{U}}_{i}^{(2)}=\frac{3}{4}\Delta
x_{i}^{(0)}\overline{\mathbf{U}}_{i}^{(0)}+\frac{1}{4}\left[\Delta
x_{i}^{(1)}\overline{\mathbf{U}}_{i}^{(1)} +\Delta
t^{n}\mathcal{L}\left(\overline{\mathbf{U}}_{i}^{(1)},
(\overline{s_{xx}})_{i}^{(1)},x_{i+\frac{1}{2}}^{(1)},w_{i+\frac{1}{2}}\right)\right],  $$
$$ (\overline{\widehat{s}_{xx}})_{i}^{(2)}=\frac{3}{4}(\overline{s_{xx}})_{i}^{(0)}
+\frac{1}{4}\left[(\overline{s_{xx}})_{i}^{(1)} +\Delta t^{n}\Theta
\left({u}_{i+\frac{1}{2}}^{(1)},(\overline{s_{xx}})_{i}^{(1)},x_{i+\frac{1}{2}}^{(1)},w_{i+\frac{1}{2}}\right)\right],  $$
$$ (\overline{s_{xx}})_{i}^{(2)}=\Upsilon((\overline{\widehat{s_{xx}}})_{i}^{(2)});  $$
\item[] {\bf Stage 3.}
$$
x_{i+\frac{1}{2}}^{(3)}=\frac{1}{3}x_{i+\frac{1}{2}}^{(0)}+\frac{2}{3}\left[x_{i+\frac{1}{2}}^{(2)}+\Delta
t^{n} \left(u_{i+\frac{1}{2}}^{(2)}-w_{i+\frac{1}{2}}\right) \right],  $$
$$  \Delta x_{i}^{(3)}=x_{i+\frac{1}{2}}^{(3)}-x_{i-\frac{1}{2}}^{(3)},  $$
$$  \Delta x_{i}^{(3)}\overline{\mathbf{U}}_{i}^{(3)}=\frac{1}{3}\Delta
x_{i}^{(0)}\overline{\mathbf{U}}_{i}^{(0)} +\frac{2}{3}\left[\Delta
x_{i}^{(2)}\overline{\mathbf{U}}_{i}^{(2)} +\Delta t^{n}\mathcal{L}\left(\overline{\mathbf{U}}_{i}^{(2)},
(\overline{s_{xx}})_{i}^{(2)},x_{i+\frac{1}{2}}^{(2)},w_{i+\frac{1}{2}}\right)\right],   $$
$$  (\overline{\widehat{s}_{xx}})_{i}^{(3)}=\frac{1}{3}(\overline{s_{xx}})_{i}^{(0)}
+\frac{2}{3}\left[(\overline{s_{xx}})_{i}^{(2)}+\Delta t^{n}\Theta\left({u}_{i+\frac{1}{2}}^{(2)},(\overline{s_{xx}})_{i}^{(2)},
x_{i+\frac{1}{2}}^{(2)},w_{i+\frac{1}{2}}\right)\right],  $$
$$ (\overline{s_{xx}})_{i}^{(3)}=\Upsilon((\overline{\widehat{s}_{xx}})_{i}^{(3)}) .$$
\end{itemize}
Here, $\mathcal{L}$ and $\Theta$ are the spatial operators representing the right-hand sides of
(\ref{e92d}) and (\ref{e108}), respectively, and
\[
x_{i+\frac{1}{2}}^{(0)}=x_{i+\frac{1}{2}}^{n},\quad
x_{i+\frac{1}{2}}^{n+1}=x_{i+\frac{1}{2}}^{(3)},\quad
u_{i+\frac{1}{2}}^{(0)}=u_{i+\frac{1}{2}}^{n},\quad
u_{i+\frac{1}{2}}^{n+1}=u_{i+\frac{1}{2}}^{(3)},
\]
\[
\Delta x_{i}^{(0)}=\Delta x_{i}^{n}, \quad
\Delta x_{i}^{n+1}=\Delta x_{i}^{(3)},
\overline{\mathbf{U}}_{i}^{(0)}=\overline{\mathbf{U}}_{i}^{n},\quad
\overline{\mathbf{U}}_{i}^{n+1}=\overline{\mathbf{U}}_{i}^{(3)},\quad
\]
\[
(\overline{s_{xx}})_{i}^{(0)}=(\overline{s_{xx}})_{i}^{n},\quad
(\overline{s_{xx}})_{i}^{n+1}=(\overline{s_{xx}})_{i}^{(3)},
\]
\begin{equation}   \label{e64}
\Upsilon(\beta)=\left\{\begin{array}{ccc}
                \beta, & for & |\beta|\leq \frac{2}{3}Y_{0} \\
                 \frac{2}{3}Y_{0}, & for & \beta>\frac{2}{3}Y_{0} \\
                -\frac{2}{3}Y_{0}, & for & \beta<-\frac{2}{3}Y_{0} .
              \end{array}\right.
\end{equation}

Recall from Section~\ref{SEC:overall} that for the free boundary points, we have
$w_{\frac{1}{2}} = w_{N+\frac{1}{2}} = 0$ by construction.
Then, at these points the above scheme reduces to
\begin{align*}
& x_{i+\frac{1}{2}}^{(1)}=x_{i+\frac{1}{2}}^{(0)}+\Delta t^{n} u_{i+\frac{1}{2}}^{n},
\qquad i = 0,\; N
\\
& x_{i+\frac{1}{2}}^{(2)}=\frac{3}{4}x_{i+\frac{1}{2}}^{(0)}+\frac{1}{4}\left[x_{i+\frac{1}{2}}^{(1)}
+\Delta t^{n} u_{i+\frac{1}{2}}^{(1)}\right],
\qquad i = 0,\; N
\\
&
x_{i+\frac{1}{2}}^{(3)}=\frac{1}{3}x_{i+\frac{1}{2}}^{(0)}+\frac{2}{3}\left[x_{i+\frac{1}{2}}^{(2)}+\Delta
t^{n} u_{i+\frac{1}{2}}^{(2)} \right],
\qquad i = 0,\; N
\end{align*}
which is exactly the third-order Runge-Kutta scheme applied to the Lagrangian movement $\dot{x} = u$.

\subsection{Time step size}
\label{SEC:CFL}

The time step size $\Delta t^{n}$ is determined using the CFL condition
$$
\Delta t^n  = 0.45 \min_{i = 1, \ldots ,N} \left({\frac{{\Delta x_i^n }}{{c_i^n+|w_i|}}} \right),
$$
where $c_i^n$ is the sound speed of solid in cell $I_i$ at
$t^n$, and $w_{i}=\frac{1}{2}(w_{i+\frac{1}{2}}+w_{i-\frac{1}{2}})$.

\section{Numerical examples}
\label{SEC:numerics}

In this section we present three numerical examples to verify the convergence order of
the MMCC scheme described in the previous section and its ability to capture
elastic-plastic waves.

\subsection{Accuracy test}

We use smooth solutions to verify the accuracy of the MMCC scheme.
The initial condition is chosen as
\[
\rho(x,0) = \rho_{0} - 0.1 \sin(2\pi x), \quad u(x,0) = 1 - 0.01 \sin(2\pi x), \quad p = 2, \quad s_{xx}=0 .
\]
The initial domain is $(0,1)$, and both ends of the domain move at the flow velocity.
A periodic boundary condition is used.
The EOS is given by the Mie-Gr\"{u}neisen model with
parameters $\rho_{0}=8930kg/m^3$, $a_{0}=3940m/s$, $\Gamma_{0}=2$,
and $s=1.49$. The constitutive model is characterized by
the shear module $\mu = 4.5\times 10^9 Pa$ and the yield
strength $Y^0=90\times 10^9 Pa$. There is no exact solution for this problem.
The reference solution is obtained using the third order cell-centered Lagrangian
scheme \cite{trrse} with a uniform mesh of $4000$ points.

The MMCC scheme is applied to the problem up to $t=1$.
The $L_1$- and $L_2$-norms of the error at the final time $t=1$ are shown in Tables \ref{msm1} and \ref{msm2}.
The results confirm the anticipated third-order accuracy of the MMCC scheme.

\begin{table}[htp]
\caption{The $L_1$-norm of the error for the MMCC scheme for the smooth solutions.} \label{msm1}
\begin{center}
\begin{tabular}{|l|c|c|c|c|c|c|c|c|}
\hline
  $N$          &   $\rho$    &  Order    &   $\rho u$  & Order     &  $\rho E$  & Order   & $s_{xx}$  & Order        \\
  \hline
    50         &  1.229E-03 &   -         &   1.680E-03 &  -       &  1.587E-03 &  -      & 7.685E-05 &  -         \\
    \hline
  100          &  4.408E-04 &  1.48       &   4.260E-04 &  1.98    &  3.451E-04 & 2.20    & 2.540E-05 & 1.60       \\
  \hline
  200          &  7.900E-05 &  2.48       &   7.430E-05 &  2.52    &  5.784E-05 & 2.58    & 4.241E-06 & 2.58       \\
  \hline
  400          &  1.054E-05 &  2.91       &   9.737E-06 &  2.93    &  6.914E-06 & 3.06    & 4.366E-07 & 3.28       \\
  \hline
\end{tabular}
\end{center}
\end{table}
%

\begin{table}[ht]
\caption{The $L_{2}$-norm of the error for the MMCC scheme for the smooth
solutions.}
\label{msm2}
\begin{center}
\begin{tabular}{|l|c|c|c|c|c|c|c|c|}
\hline
 $N$          &   $\rho$   &  Order       &   $\rho u$  & Order   &  $\rho E$  & Order   & $s_{xx}$  & Order         \\
 \hline
  50           &  1.414E-03 &   -          &   2.187E-03 &  -         & 2.095E-03 &  -       & 8.911E-05 &  -         \\
   \hline
  100          &  5.148E-04 &  1.46        &   5.852E-04 & 1.90       & 5.104E-04 & 2.04     & 3.276E-05 & 1.44        \\
   \hline
  200          &  1.028E-04 &  2.32        &   1.164E-04 & 2.33       & 9.902E-05 & 2.37     & 6.068E-06 & 2.43        \\
   \hline
  400          &  1.534E-05 &  2.74        &   1.681E-05 & 2.79       & 1.302E-05 & 2.93     & 7.097E-07 & 3.10        \\
\hline
\end{tabular}
\end{center}
\end{table}

\subsection{A piston problem with stress shock waves in copper}

This piston problem is concerned with a piece of copper with the following initial setting:
the length of copper is $1m$, the initial density
$\rho_{0}=8930kg/m^3$, the initial pressure $p_{0}=10^5Pa$, and the
initial velocity is zero. The EOS for copper is given
by the Mie-Gr\"{u}neisen model with parameters
$\rho_{0}=8930kg/m^3$, $a_{0}=3940m/s$, $\Gamma_{0}=2$, and $s=1.49$.
The constitutive model is characterized by the shear module $\mu = 45\times 10^9 Pa$
and the yield strength $Y^0=90\times 10^6 Pa$. A velocity boundary condition $v_{piston}=20m/s$
is used on the left free boundary while the wall boundary condition is implemented
on the right fixed boundary. This problem has an exact solution; see \cite{cc3}.

\begin{figure}[htp]
\begin{center}
\includegraphics[width=10cm]{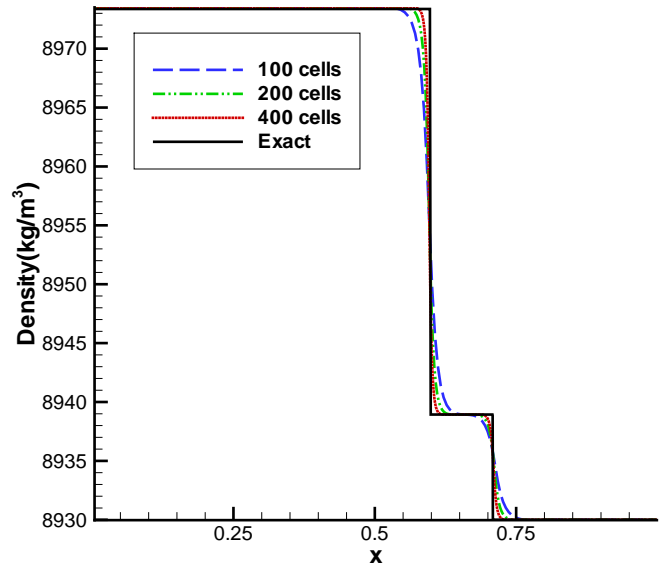}
\includegraphics[width=10cm]{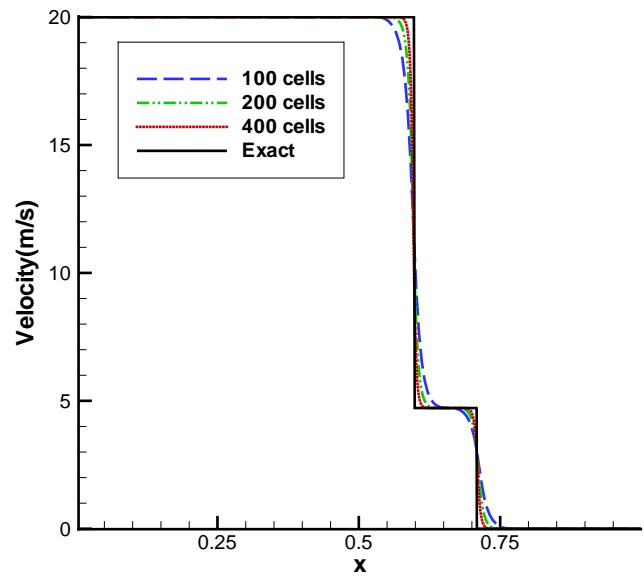}
\end{center}
\caption{Piston problem. The density and velocity at $t=150\mu s$ are shown.}    \label{fig1-dv}
\end{figure}
\begin{figure}[htp]
\begin{center}
\includegraphics[width=10cm]{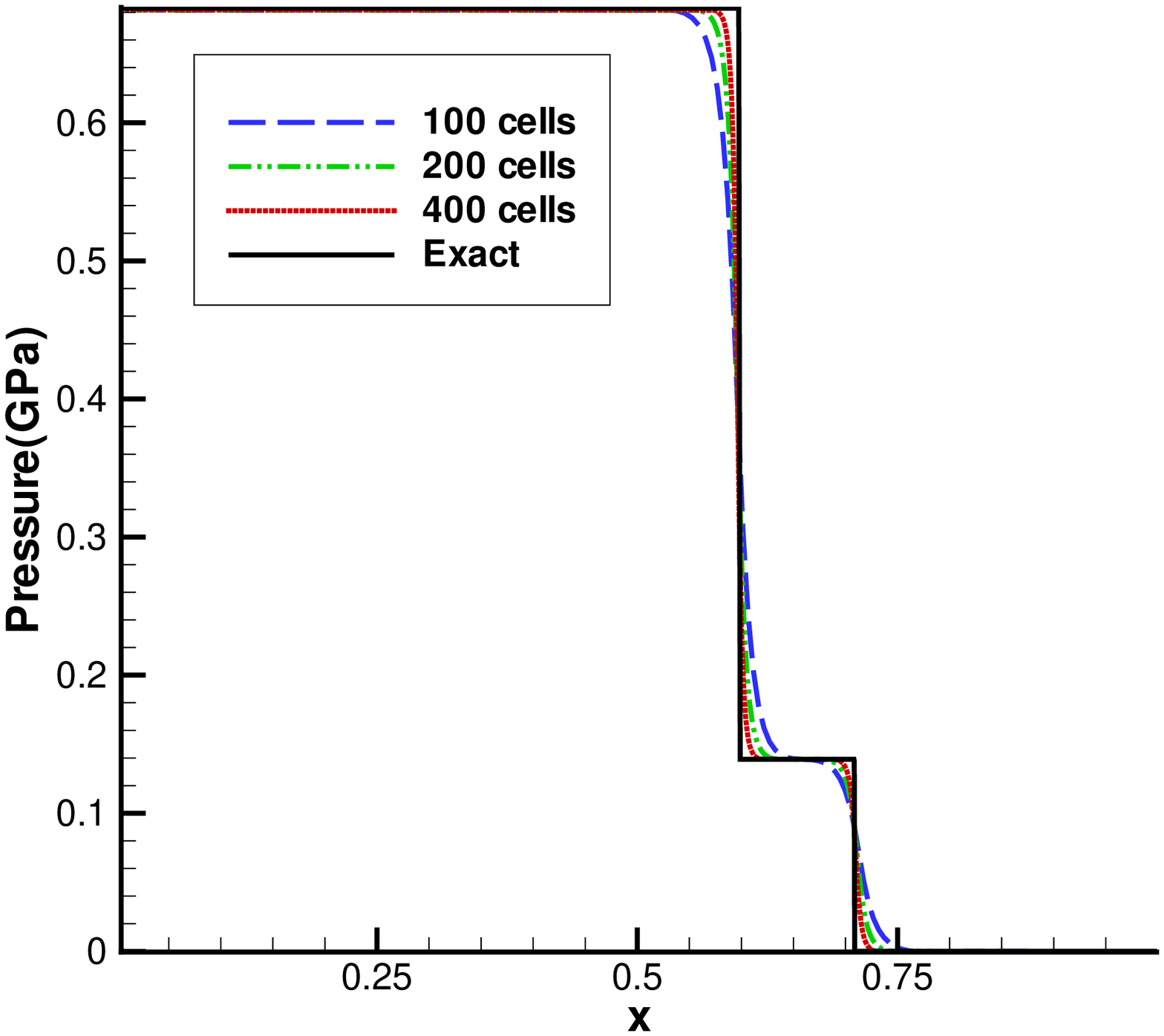}
\includegraphics[width=10cm]{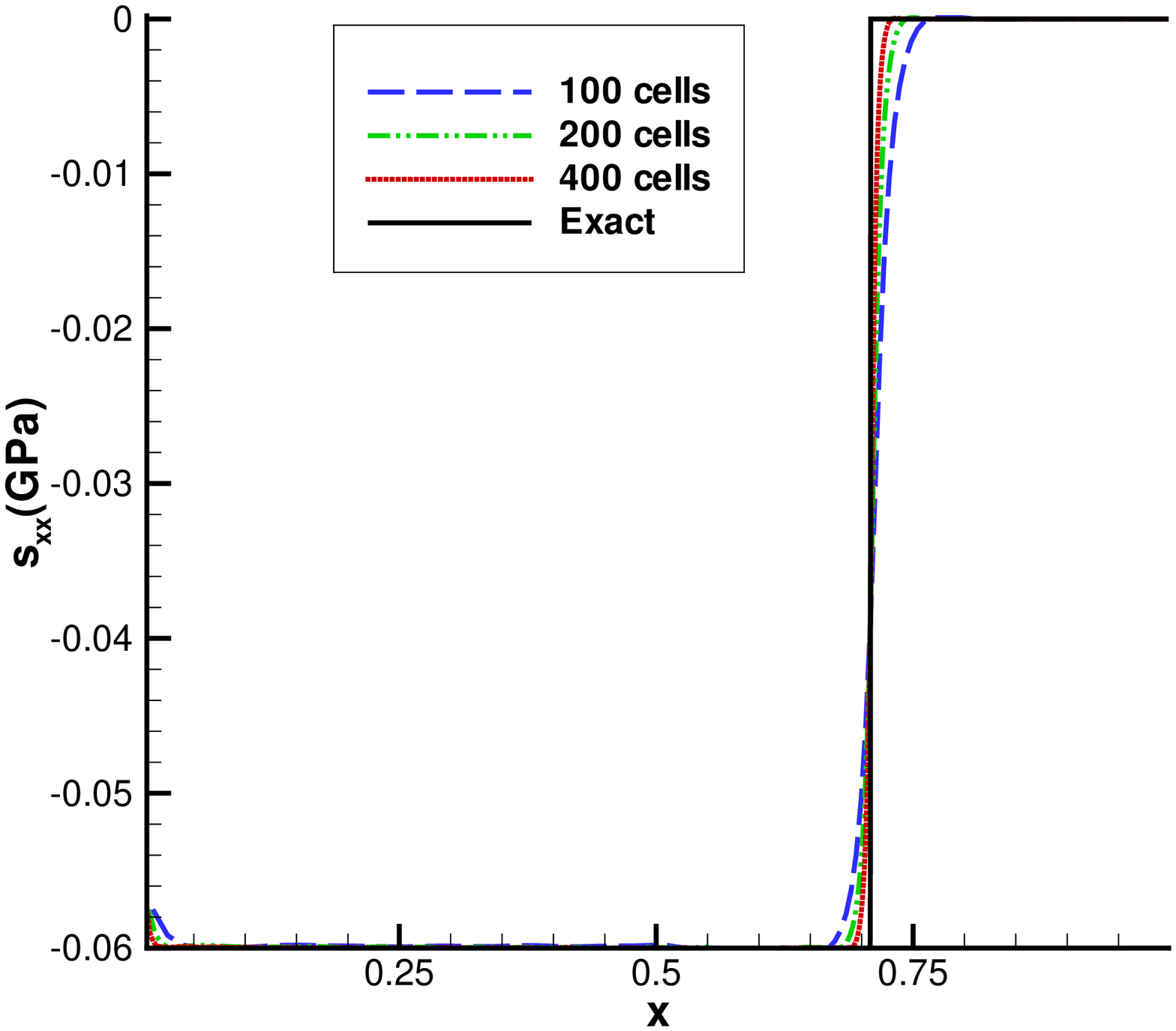}
\end{center}
\caption{Piston problem. The pressure and $s_{xx}$ at $t=150\mu s$ are shown.}   \label{fig1-ps}
\end{figure}

To show the convergence of the MMCC scheme, we solve the problem
with $100$, $200$, and $400$ cells.
The final time is $t=150\mu s$. The numerical results with
different $N$ are shown in Figs.~\ref{fig1-dv} and \ref{fig1-ps}.
From the figures we can see that the numerical solution is getting closer to
the exact solution for larger $N$.
Moreover, the scheme well captures the plastic shock wave (near $x = 0.6$) and
the elastic shock wave (near $x= 0.7$)
and there is no numerical oscillation near the shock waves.
Furthermore, the yielding condition (\ref{e11}) is in effect for $ 0 < x < 0.7$
(cf. Fig.~\ref{fig1-ps}). It creates the plastic shock wave and also makes
$s_{xx}$ to be only piecewise smooth.
The $L_1$- norm of the error at the final time
$t=150 \mu s$ are shown in Table \ref{piston}. Due to the discontinuity of the exact solution
at the elastic shock wave and its non-smoothness at the plastic shock wave, the scheme cannot
be expected to converge at the optimal third-order rate.
The results show that the convergence is about first order for this problem as the mesh is refined.

\begin{figure}[htp]
\begin{center}
\includegraphics[width=10cm]{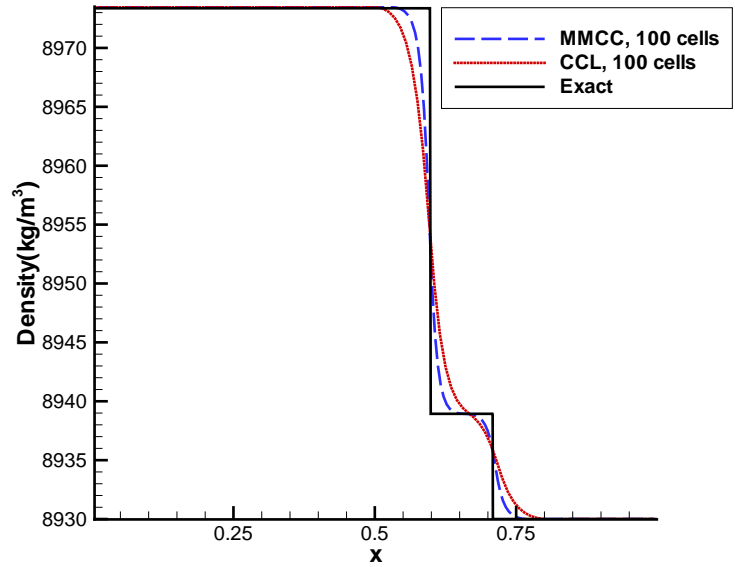}
\includegraphics[width=10cm]{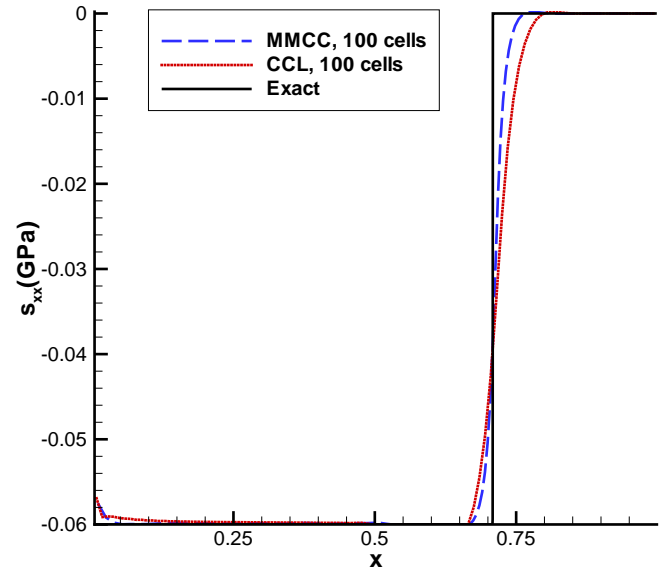}
\end{center}
\caption{Piston problem. The density and $s_{xx}$ are shown for MMCC and CCL schemes
with 100 cells.} \label{compar1-d}
\end{figure}

\begin{figure}[ht]
\begin{center}
\includegraphics[width=6cm]{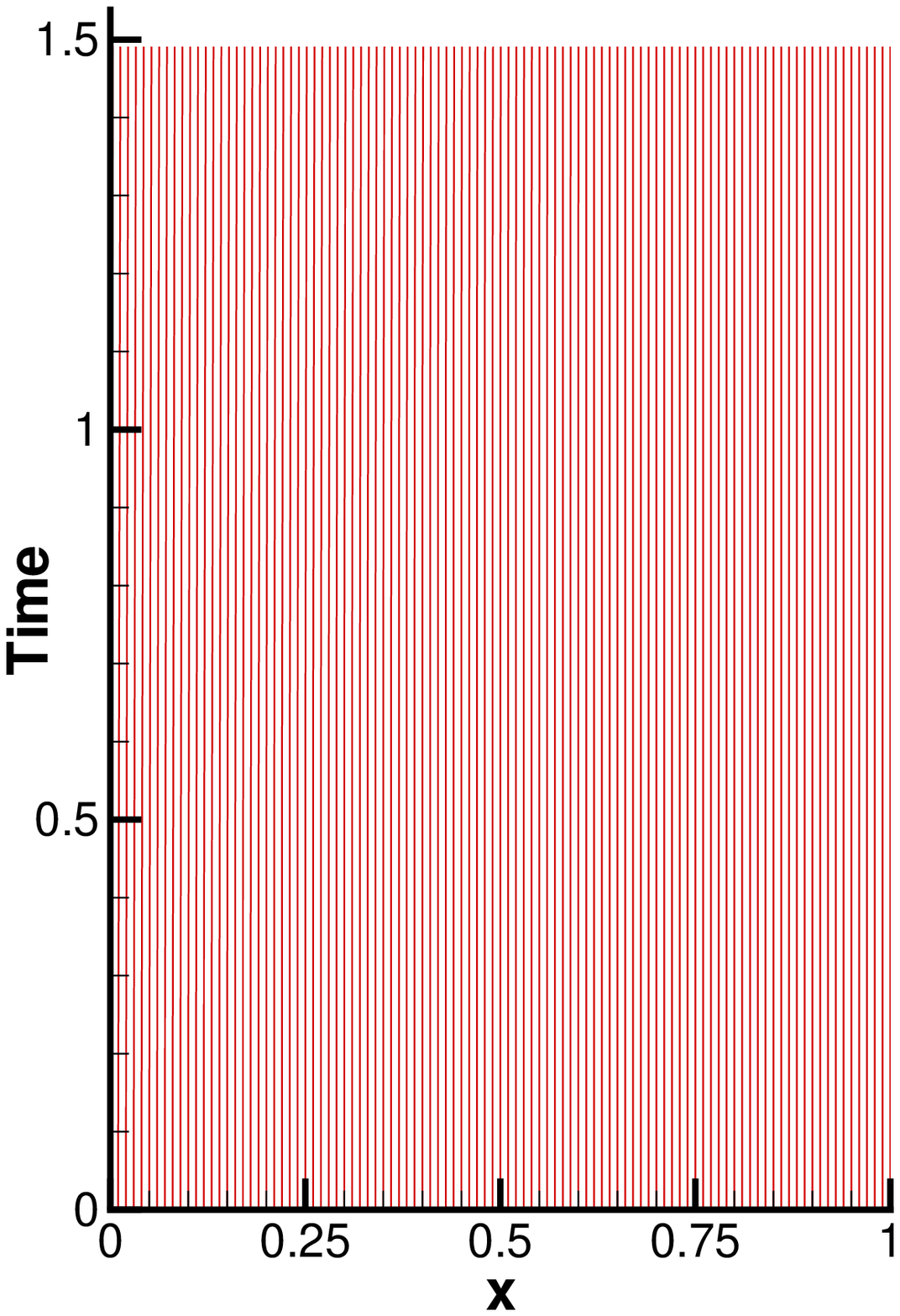}
\includegraphics[width=6cm, height = 9cm]{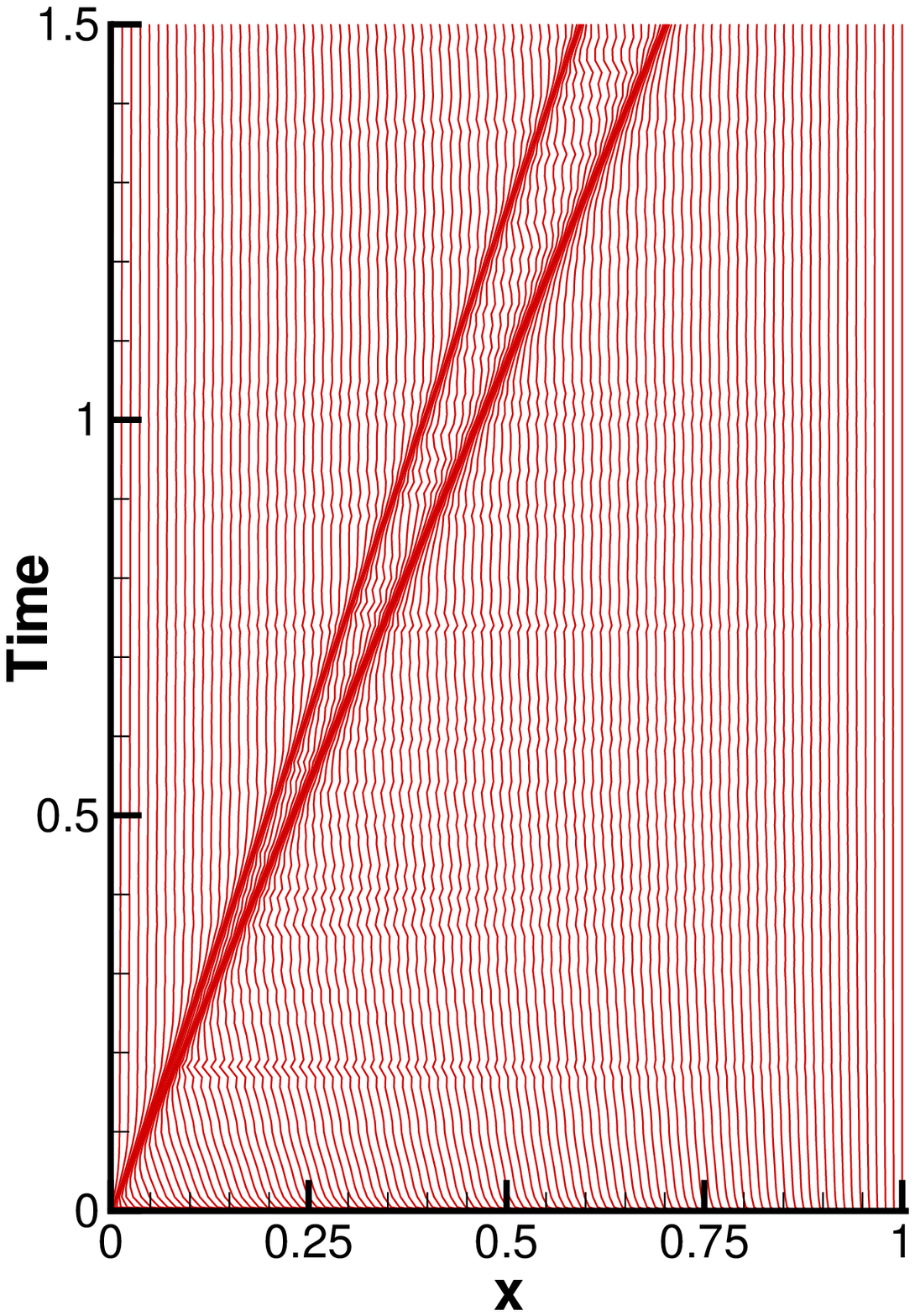}
\end{center}
\caption{Piston problem. Mesh trajectories for CCL (left) and MMCC (right) schemes with 100 cells.}    \label{fig1-xt}
\end{figure}

We also compare the MMCC scheme with the third-order cell-centered Lagrangian scheme (CCL)
based on the same Riemann solver TRRSE of \cite{trrse}. The comparison results
are shown in Fig.~\ref{compar1-d}. We can see that MMCC is more accurate than CCL especially
around the shock waves. This is due to the fact that MMCC concentrates more points around the shock waves;
see the mesh trajectories shown in Fig. \ref{fig1-xt}. No mesh crossing has been experienced
for all of the computation for this example.

\begin{table}[ht]
\caption{The $L_{1}$-norm of the error for the MMCC scheme for the piston
problem.}
\label{piston}
\begin{center}
\begin{tabular}{|l|c|c|c|c|c|c|c|c|}
\hline
  $N$          &   $\rho$   &  Order       &   $\rho u$  & Order     &  $\rho E$  & Order   & $s_{xx}$ & Order        \\
 \hline
  100          &  6.732E-04 &    -         &   2.783E-04 &   -       & 7.429E-07 &   -     & 1.038E-05 &   -       \\
\hline
  200          &  2.871E-04 &  1.23        &   1.195E-04 & 1.22      & 3.817E-07 & 0.96    & 4.924E-06 & 1.08       \\
 \hline
  400          &  1.249E-04 &  1.20        &   5.207E-05 & 1.20      & 1.856E-07 & 1.04    & 2.316E-06 & 1.09       \\
\hline
\end{tabular}
\end{center}
\end{table}


%

\subsection{Wilkins' problem with the Mie-Gr\"{u}neisen  EOS}

This problem, originally introduced by Wilkins \cite{wilkins}, is
used here to show the ability of the MMCC scheme to compute rarefaction
waves. It describes a moving aluminium plate striking on
another aluminium plate. The EOS for aluminium is the
Mie-Gr\"{u}neisen model with parameters $\rho_{0}=2785kg/m^3$,
$a_{0}=5328m/s$, $\Gamma_{0}=2$, and $s=1.338$. The constitutive model
is characterized by the shear module $\mu =
27.6\times 10^9 Pa$ and the yield strength $Y^0=300 \times 10^6 Pa$.
The initial condition of this problem is
\begin{align}
& (\rho,u,p) =
 \begin{cases}
                      (2785kg/m^3,\; 800m/s,\; 10^{-6}Pa), & \mbox{for}\;\; 0m \leq x \leq 5\times 10 ^ {-3}m, \\[1mm]
                      (2785kg/m^3,\; 0m/s,\; 10^{-6}Pa), & \mbox{for}\;\; 5\times 10 ^ {-3}m \leq x \leq  50\times 10 ^
                      {-3}m.
                    \end{cases}
\end{align}

\begin{figure}[htp]
\begin{center}
\includegraphics[width=10cm]{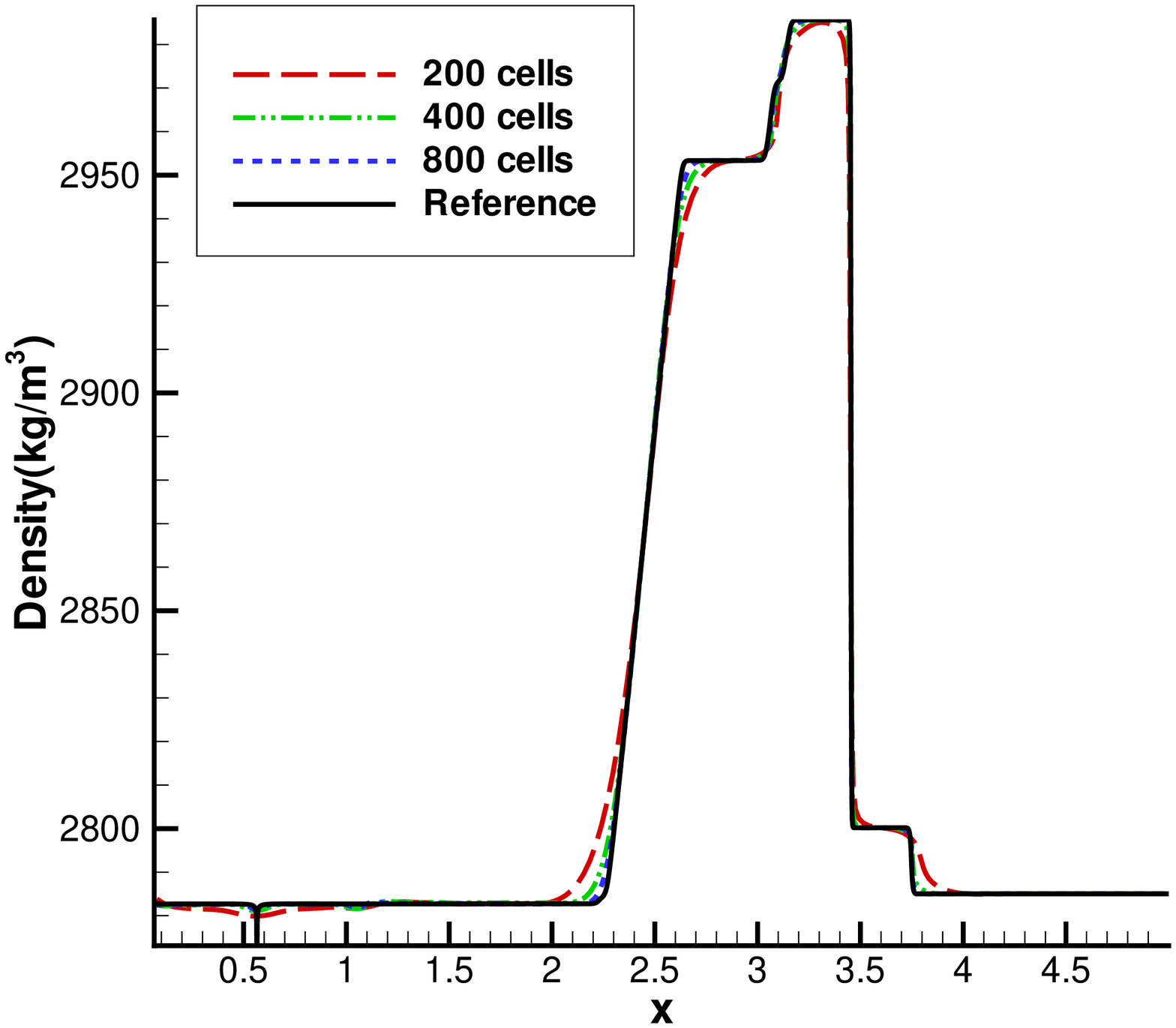}\\
\includegraphics[width=8cm,height=8cm]{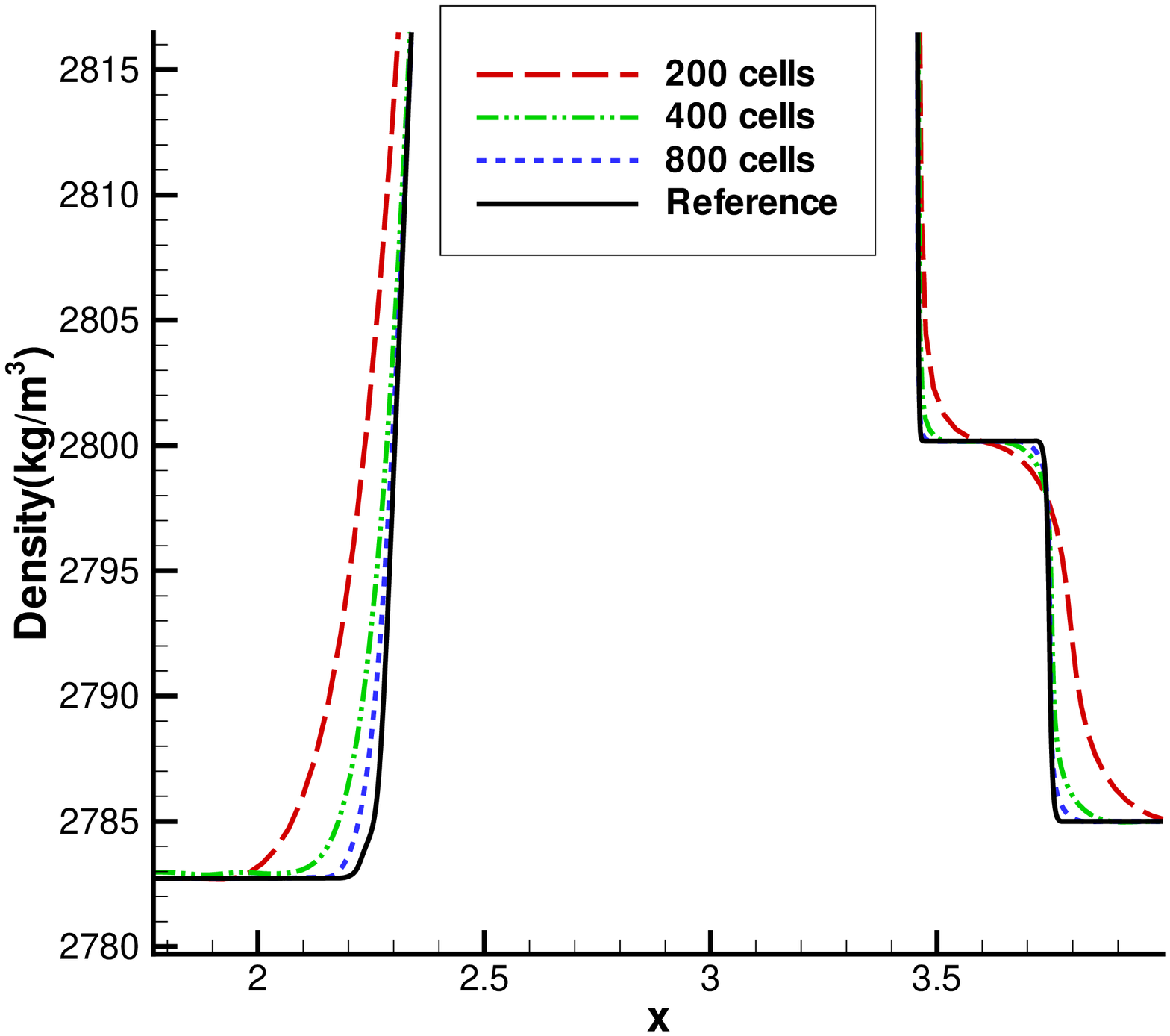}
\includegraphics[width=8cm,height=8cm]{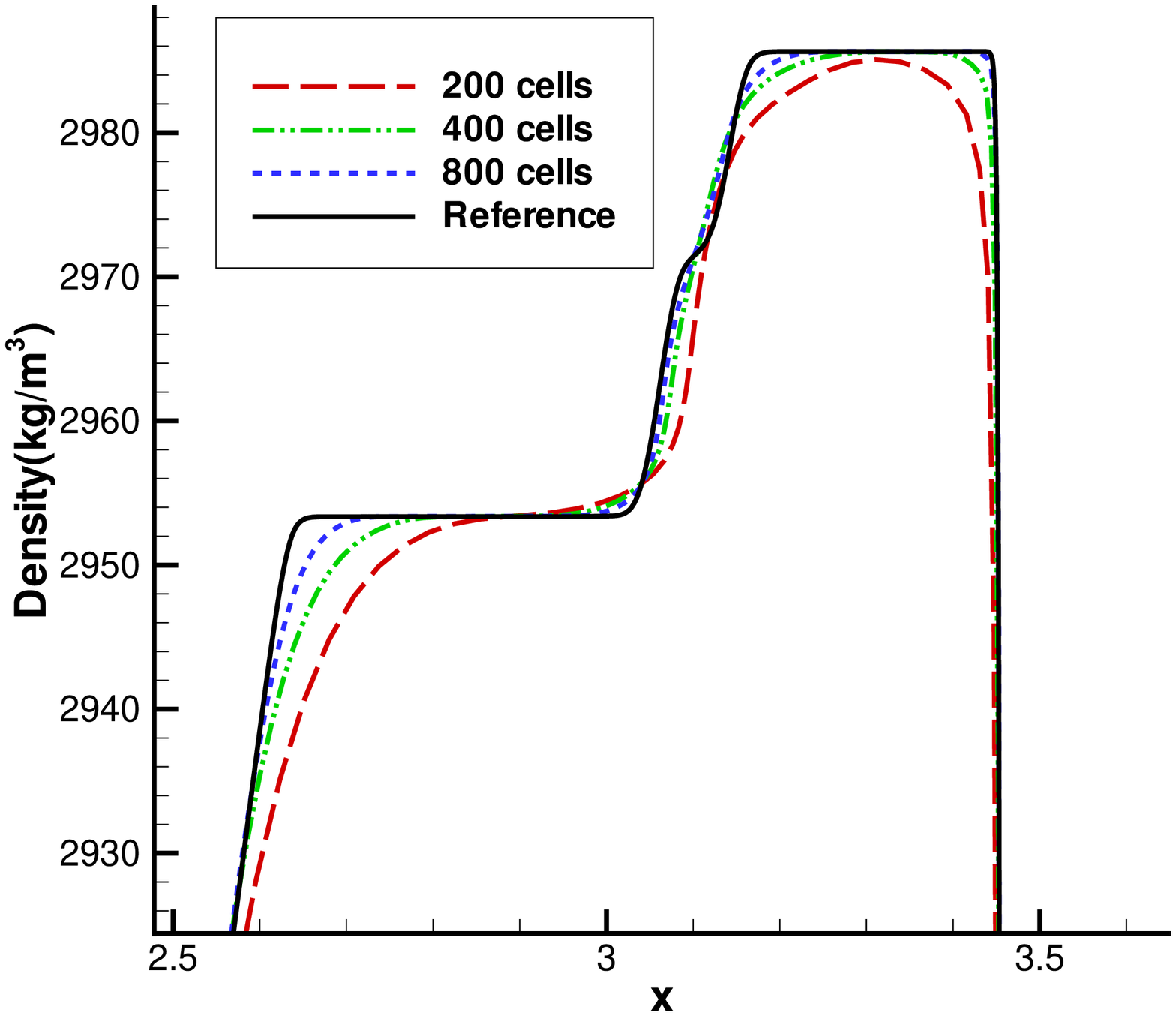}
\end{center}
\caption{Wilkins' problem with the Mie-Gr\"{u}neisen EOS.
The density at time $t=5\mu s$ is shown. The bottom figures are the close-up of
the regions around the rarefaction and shock waves.}  \label{fig2-d}
\end{figure}

\begin{figure}[htp]
\begin{center}
\includegraphics[width=10cm]{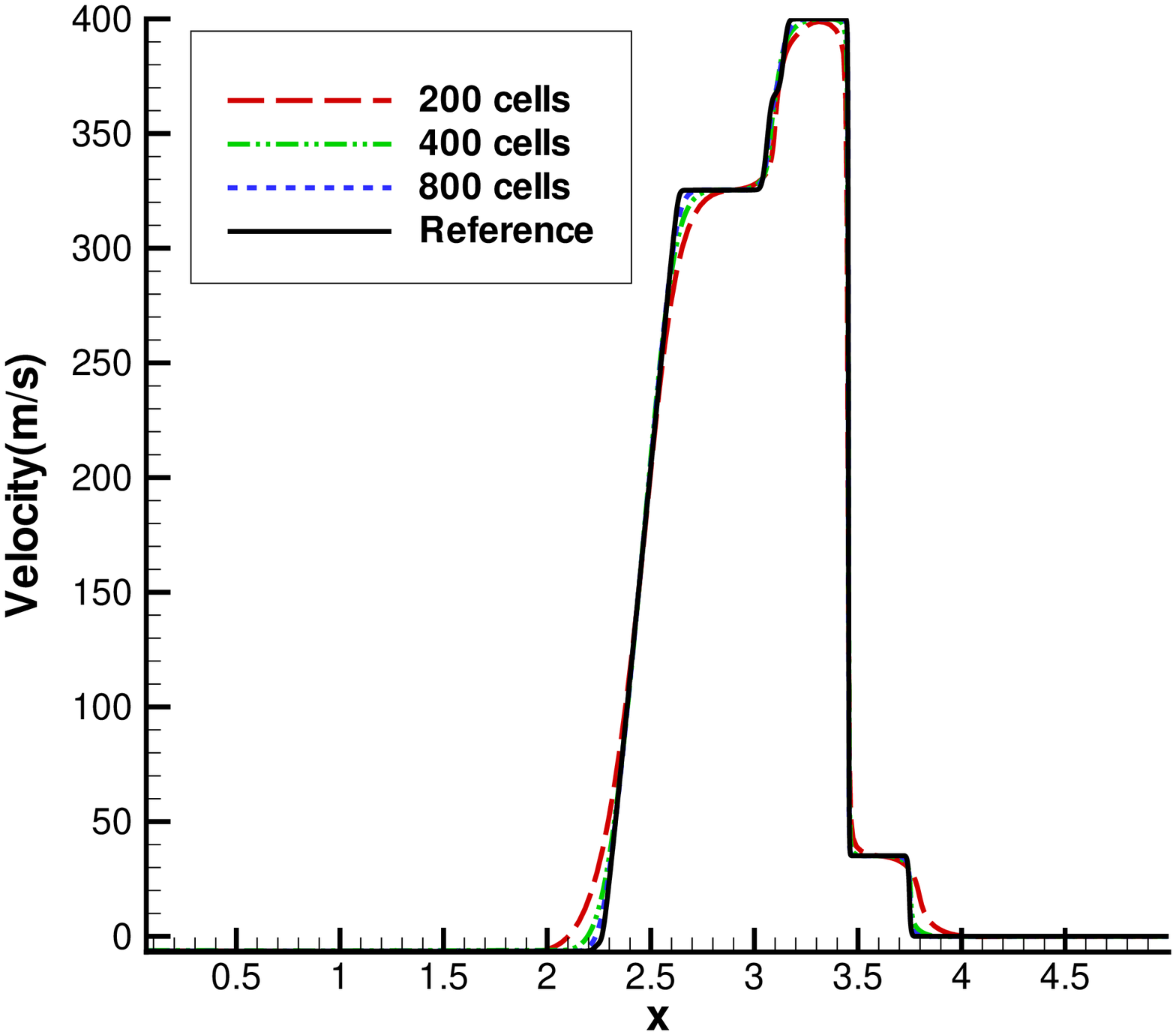}\\
\includegraphics[width=8cm,height=8cm]{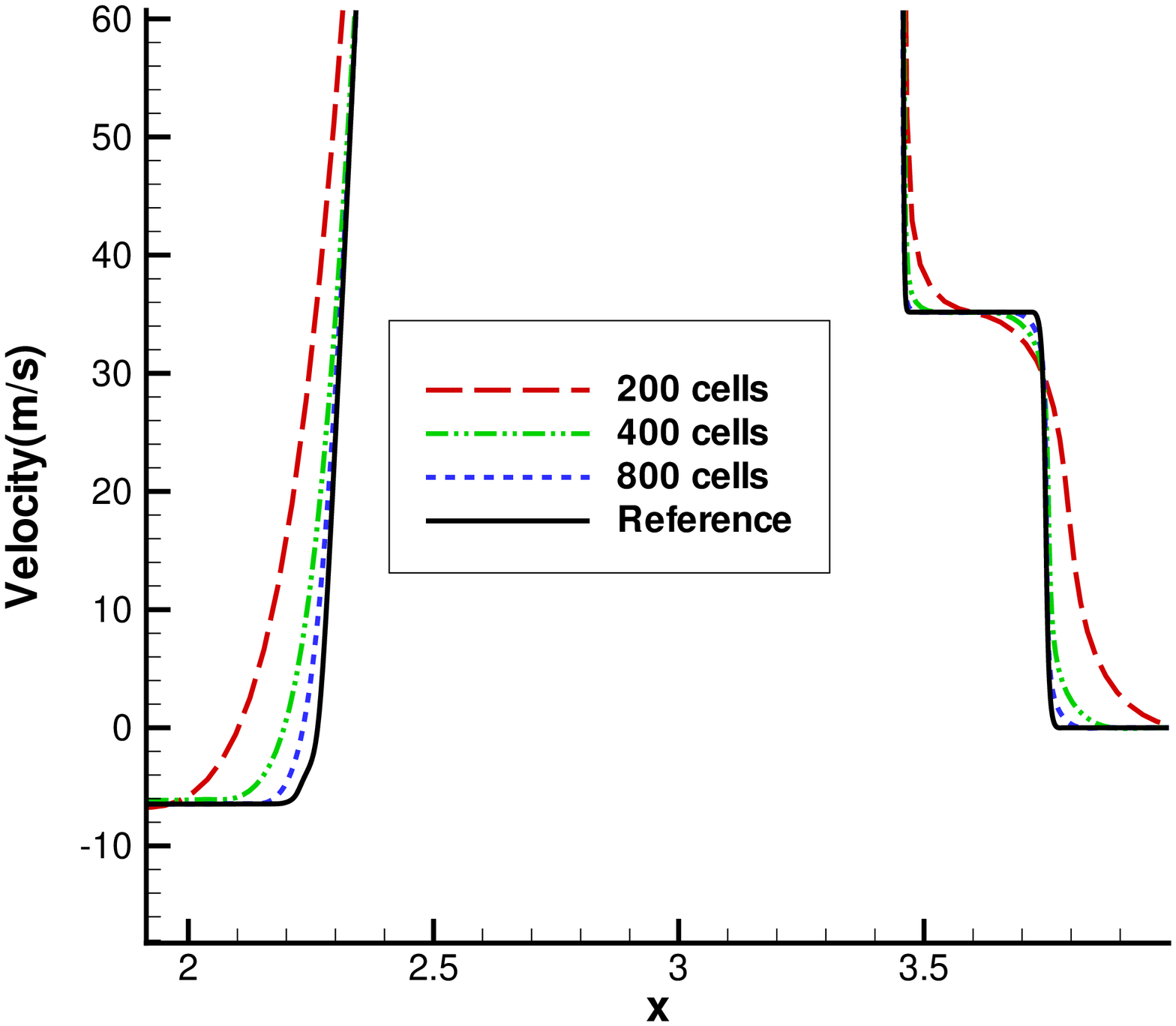}
\includegraphics[width=8cm,height=8cm]{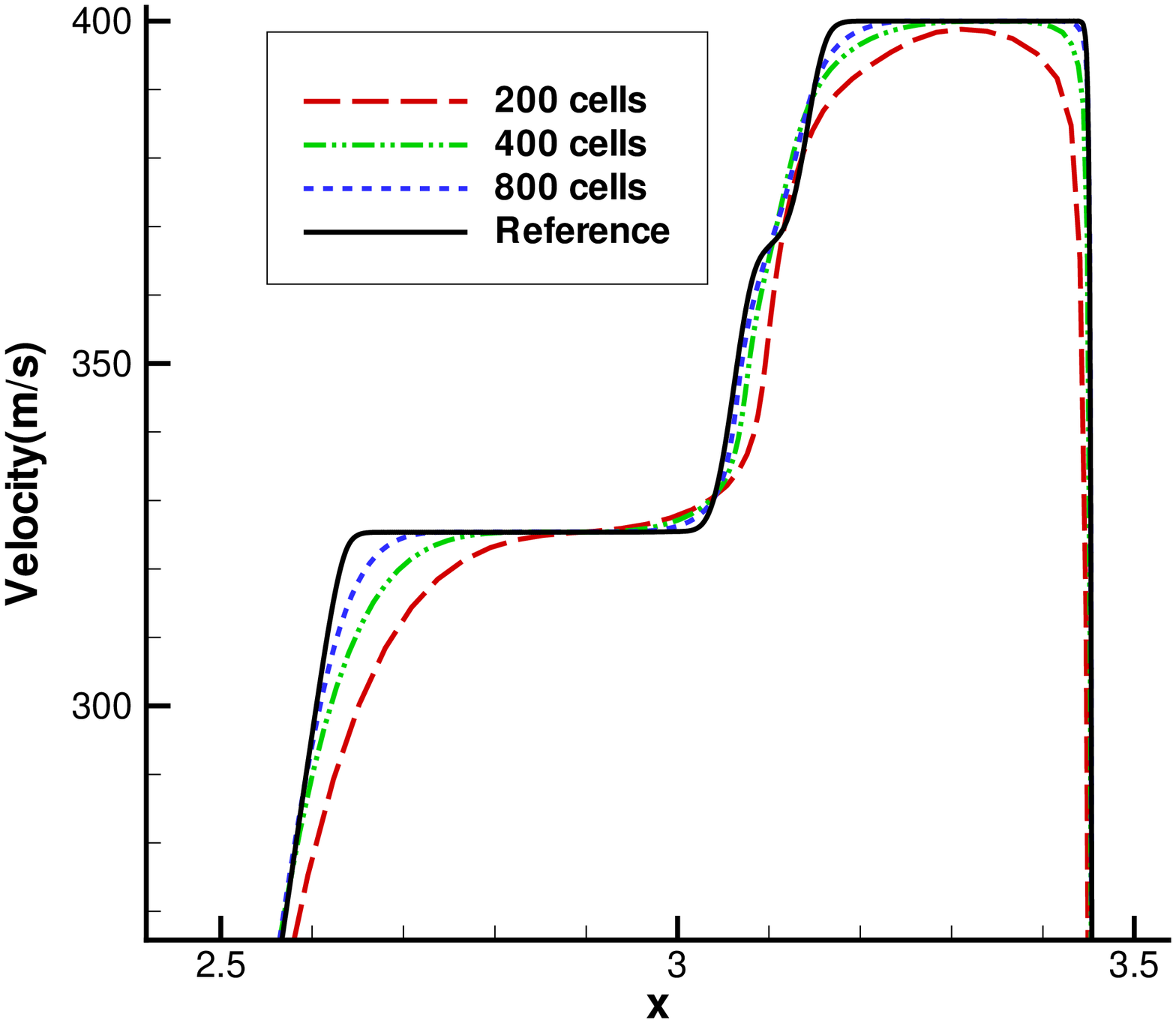}
\end{center}
\caption{Wilkins' problem with the Mie-Gr\"{u}neisen  EOS. The velocity
at $t=5\mu s$ is shown. The bottom figures are the close-up of
the regions around the rarefaction and shock waves.}  \label{fig2-v}
\end{figure}

\begin{figure}[htp]
\begin{center}
\includegraphics[width=10cm]{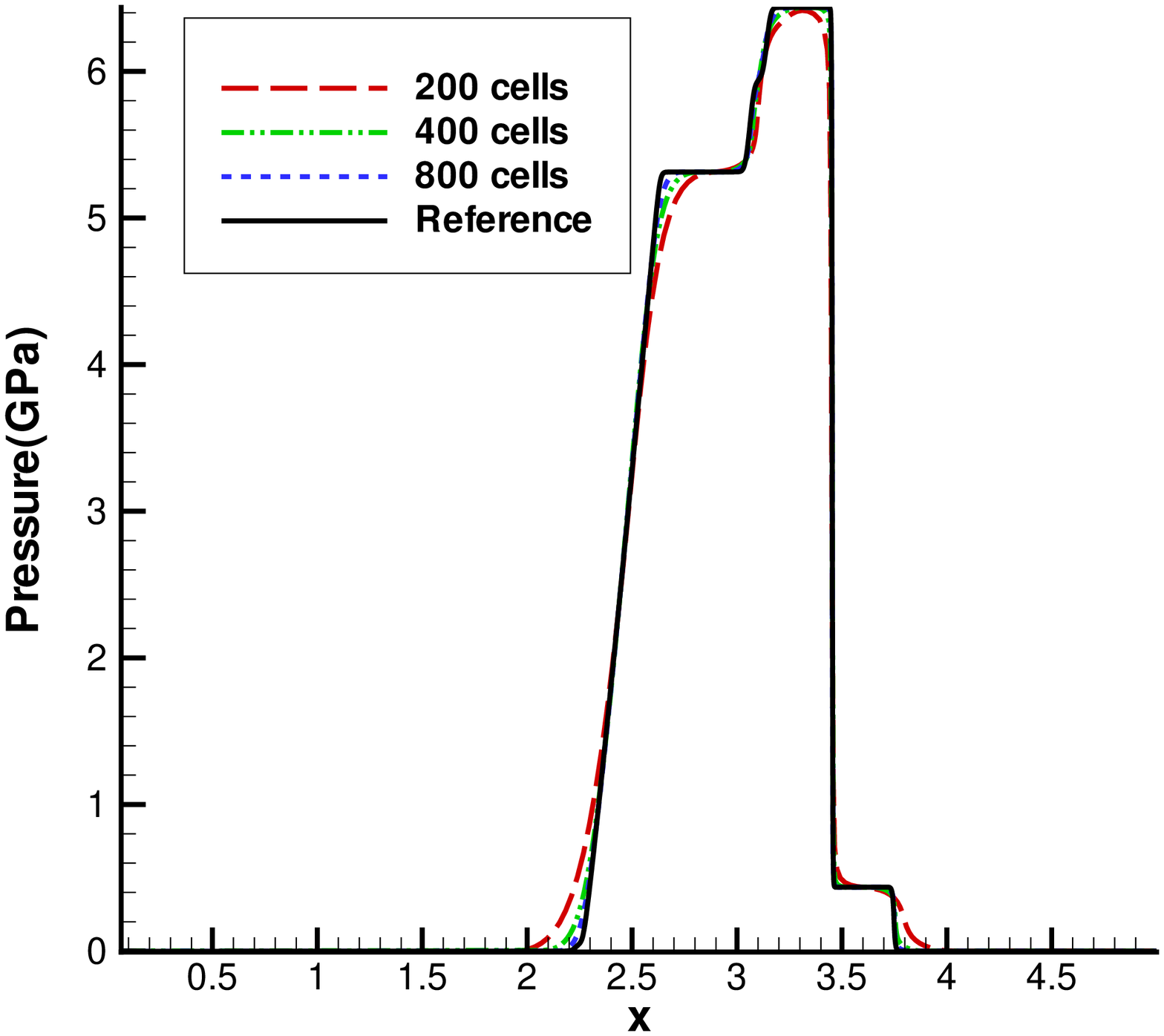}\\
\includegraphics[width=8cm,height=8cm]{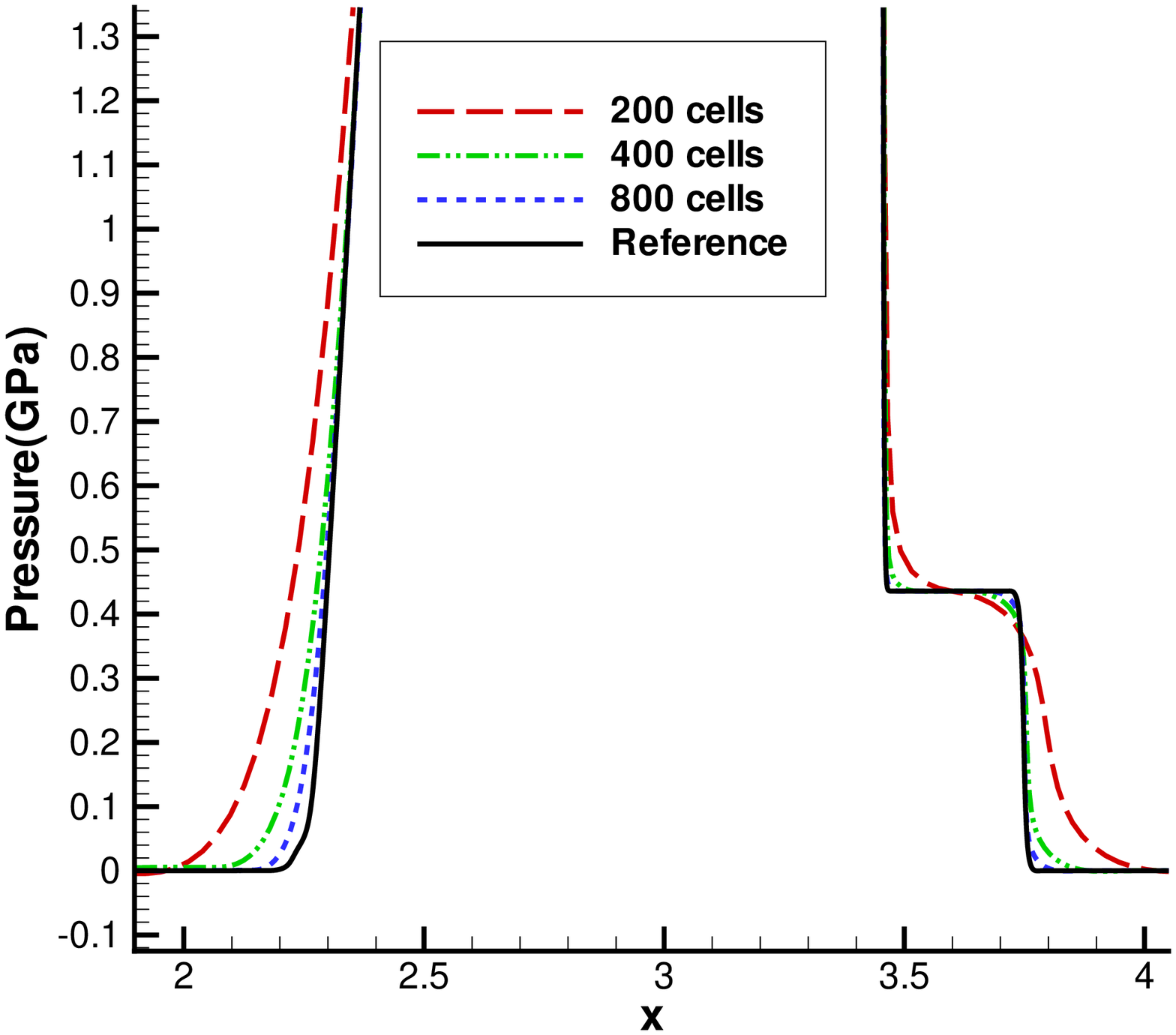}
\includegraphics[width=8cm,height=8cm]{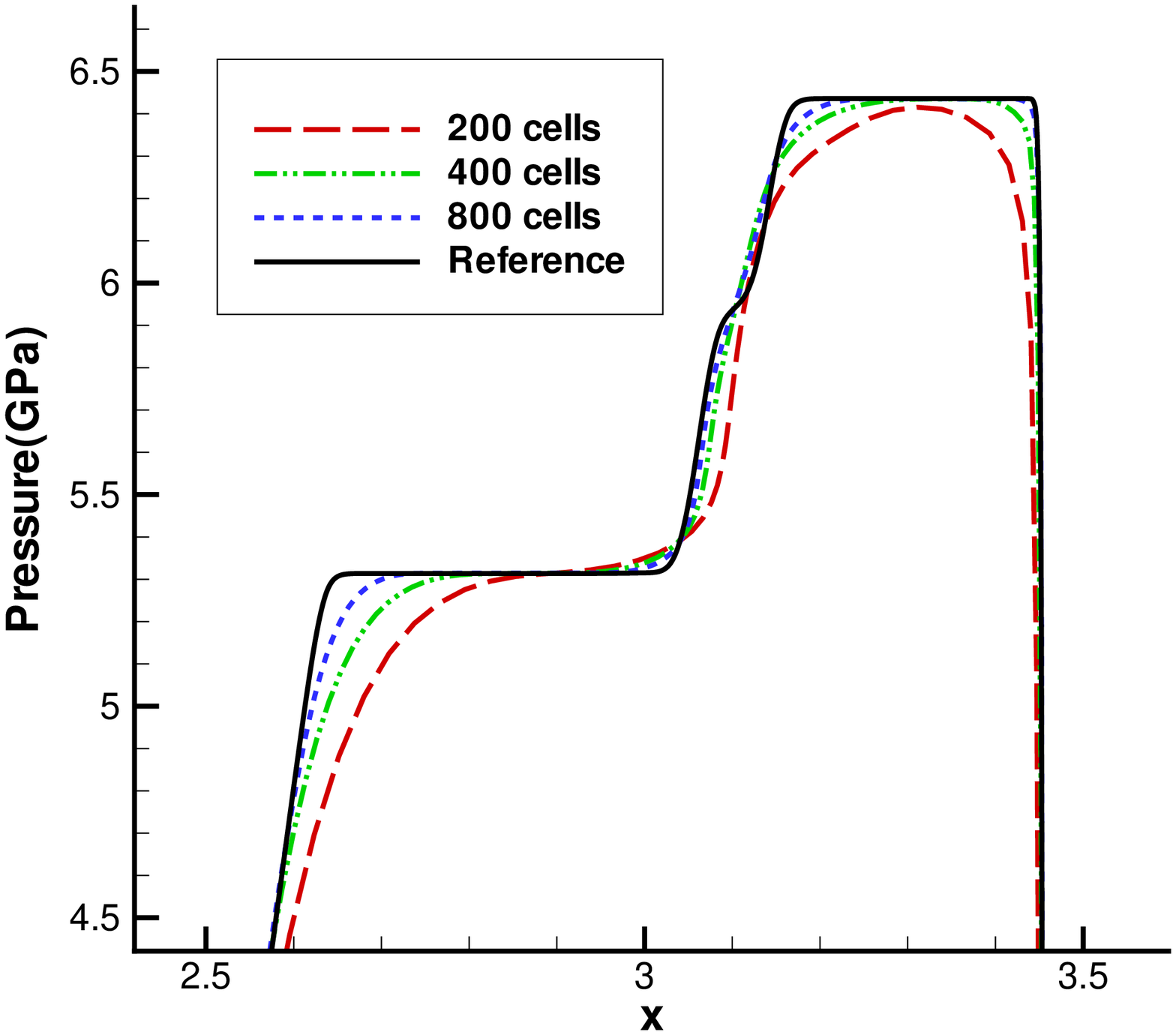}
\end{center}
\caption{Wilkins' problem with the Mie-Gr\"{u}neisen  EOS. The pressure
at $t=5\mu s$ is shown. The bottom figures are the close-up of
the regions around the rarefaction and shock waves.}  \label{fig2-p}
\end{figure}

\begin{figure}[htp]
\begin{center}
\includegraphics[width=10cm]{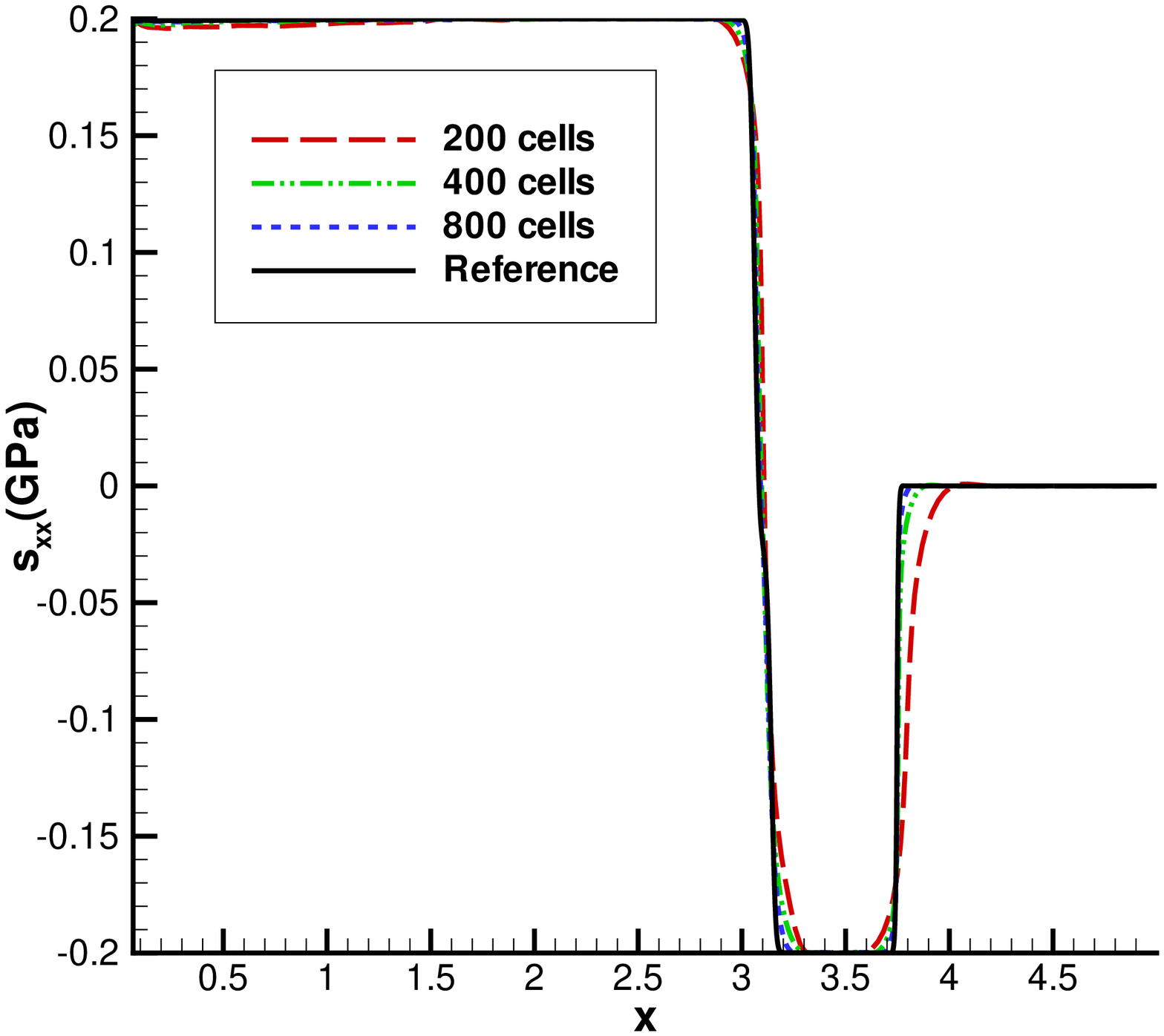}\\
\includegraphics[width=8cm]{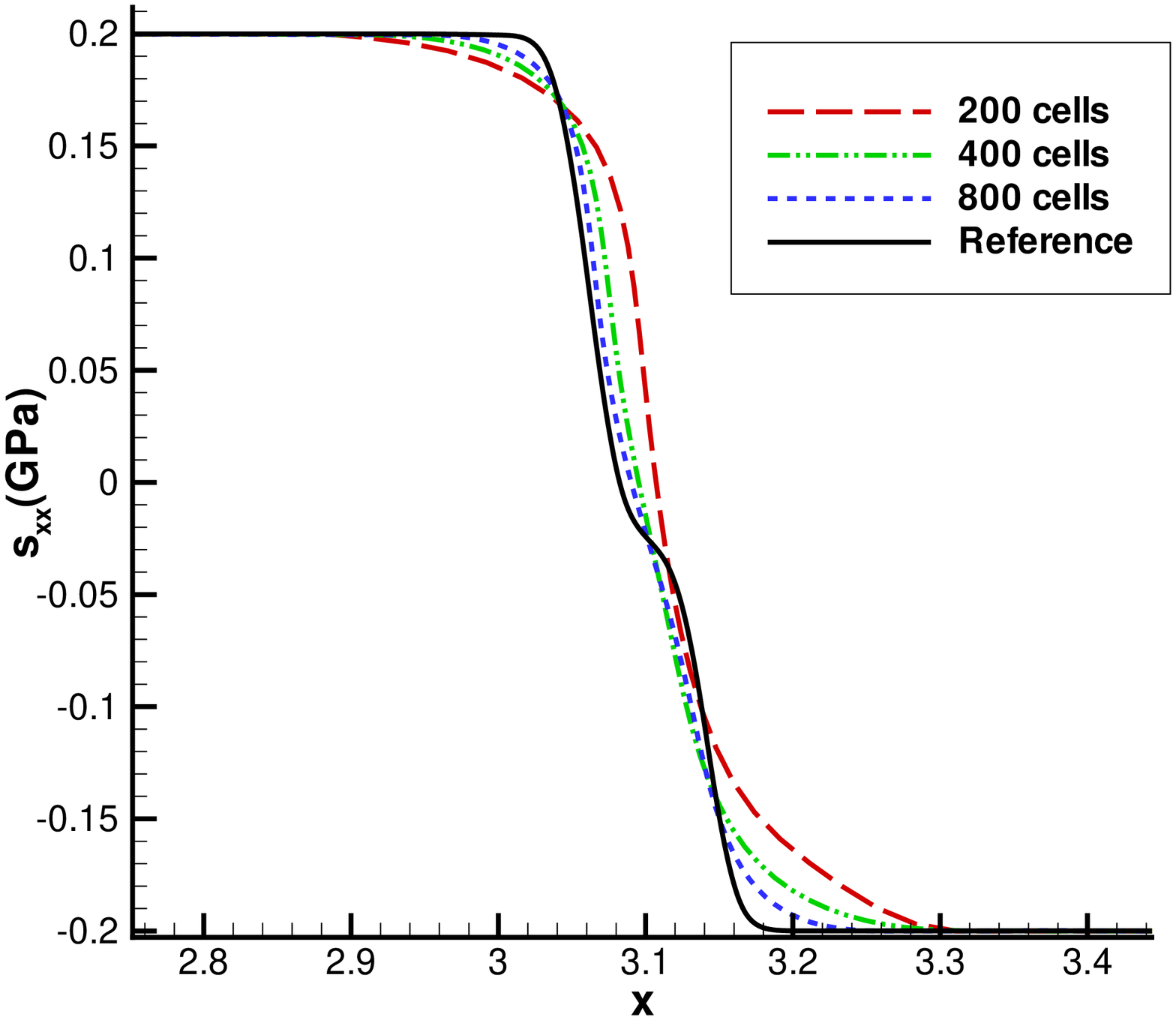}
\includegraphics[width=8cm]{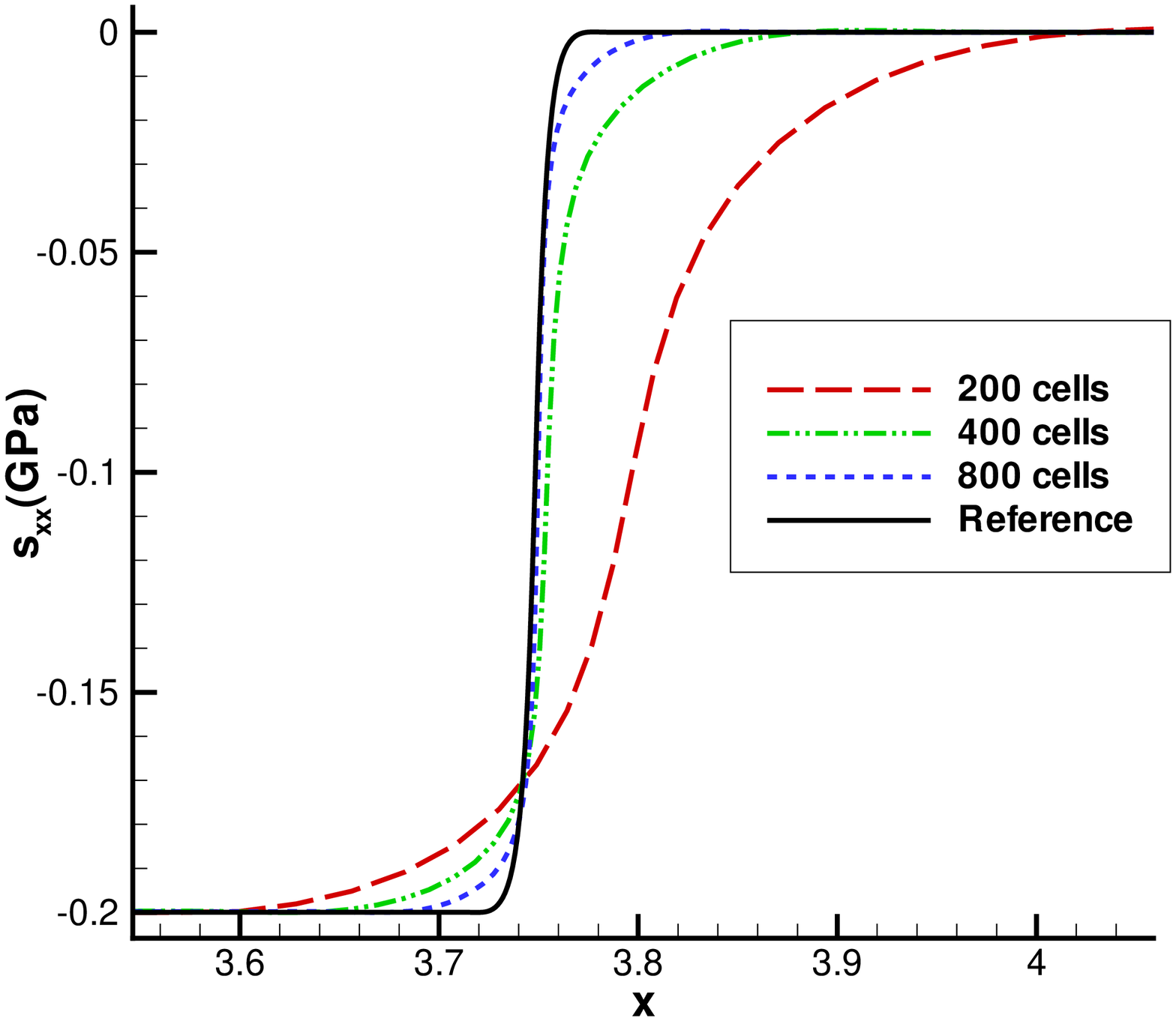}
\end{center}
\caption{Wilkins' problem with the Mie-Gr\"{u}neisen  EOS. $s_{xx}$
at $t=5\mu s$ is shown. The bottom figures are the close-up of
the region around the right-facing elastic shock wave.}
\label{fig2-s}
\end{figure}

We use $200$, $400$, and $800$ cells to solve the problem up to
$ t = 5\mu s$ with a free boundary condition on the left boundary
and the wall boundary condition on the right boundary. Computed results
are presented in Figs.~\ref{fig2-d} to \ref{fig2-s}. The reference
solution (solid lines in all figures) is computed by the third-order cell-centered Lagrangian scheme
based on HLLCE Riemann solvers given in \cite{HLLCE}
using $4000$ cells.

From Figs.~\ref{fig2-d} to \ref{fig2-s} we can observe that the
numerical solution converges to the reference one as the number of cells increases.
Moreover, the elastic and plastic right-facing shocks and the reflected elastic
and plastic rarefaction waves are resolved well. These results are in good agreement with
those in \cite{HLLCE}.

\begin{figure}[htp]
\begin{center}
\includegraphics[width=12cm]{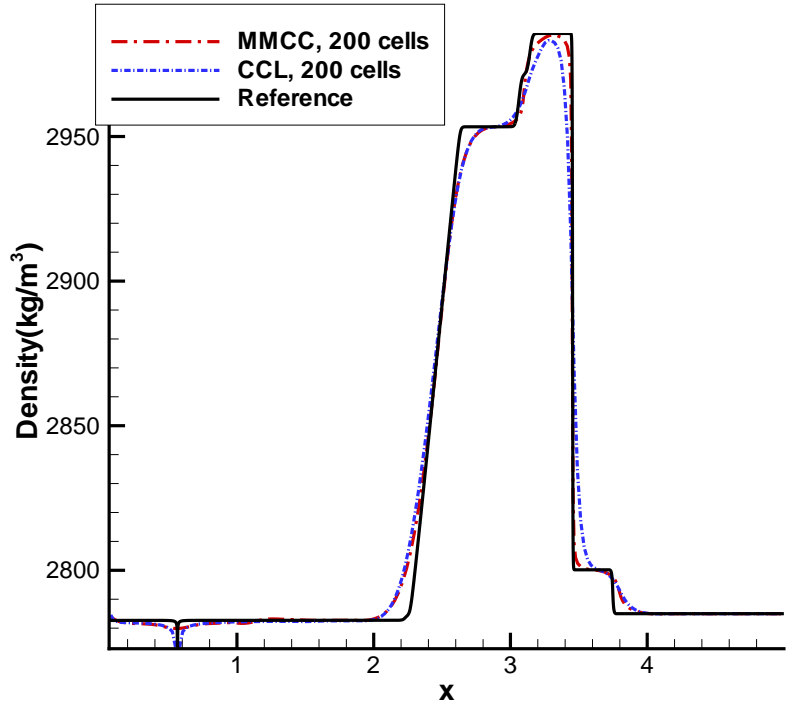}\\
\includegraphics[width=12cm]{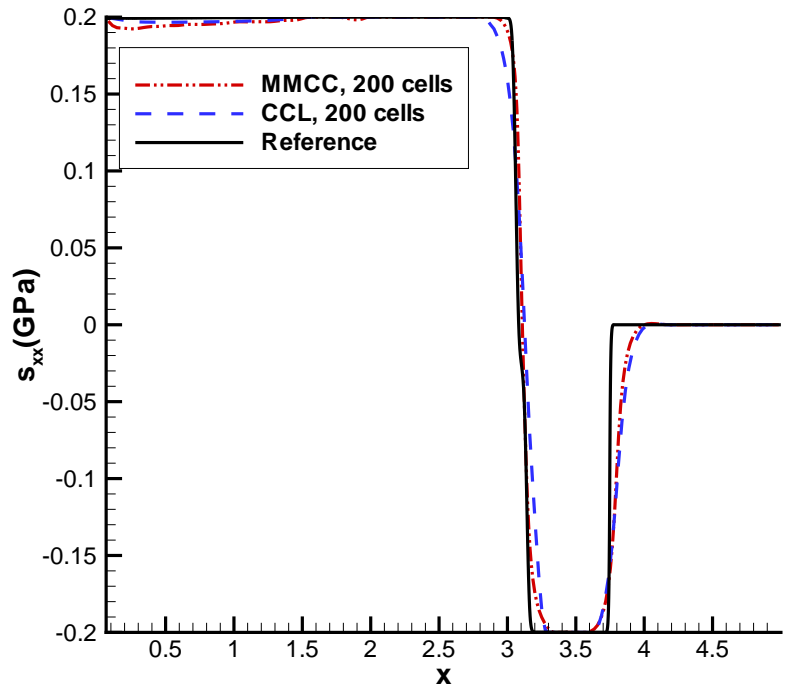}
\end{center}
\caption{Wilkins' problem with the Mie-Gr\"{u}neisen EOS. The density and $s_{xx}$ are shown for
MMCC and CCL with 200 cells. }
\label{compar2-d}
\end{figure}

\begin{figure}[ht]
\begin{center}
\includegraphics[width=6cm, height=7.8cm]{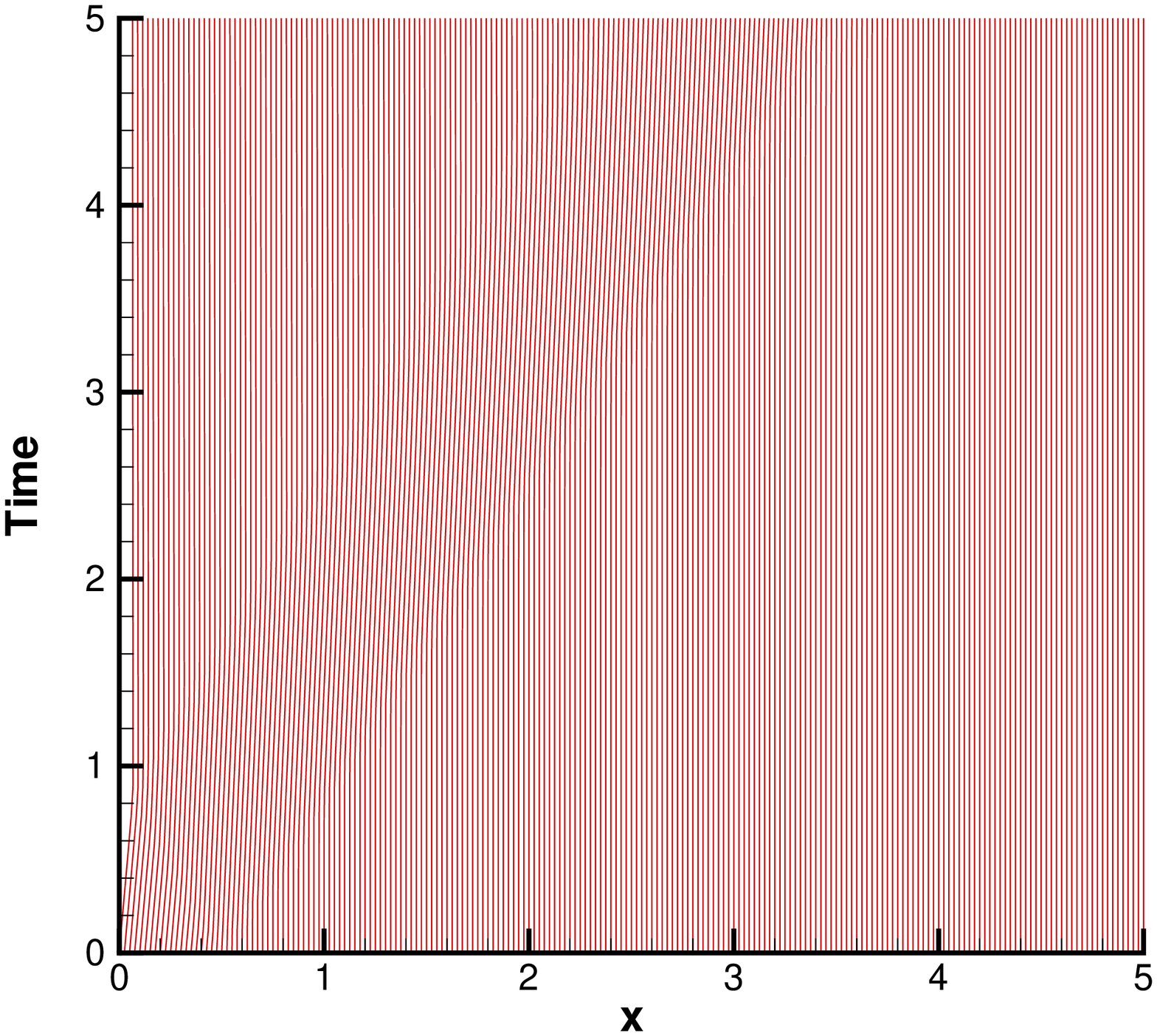}
\includegraphics[width=6cm, height=8cm]{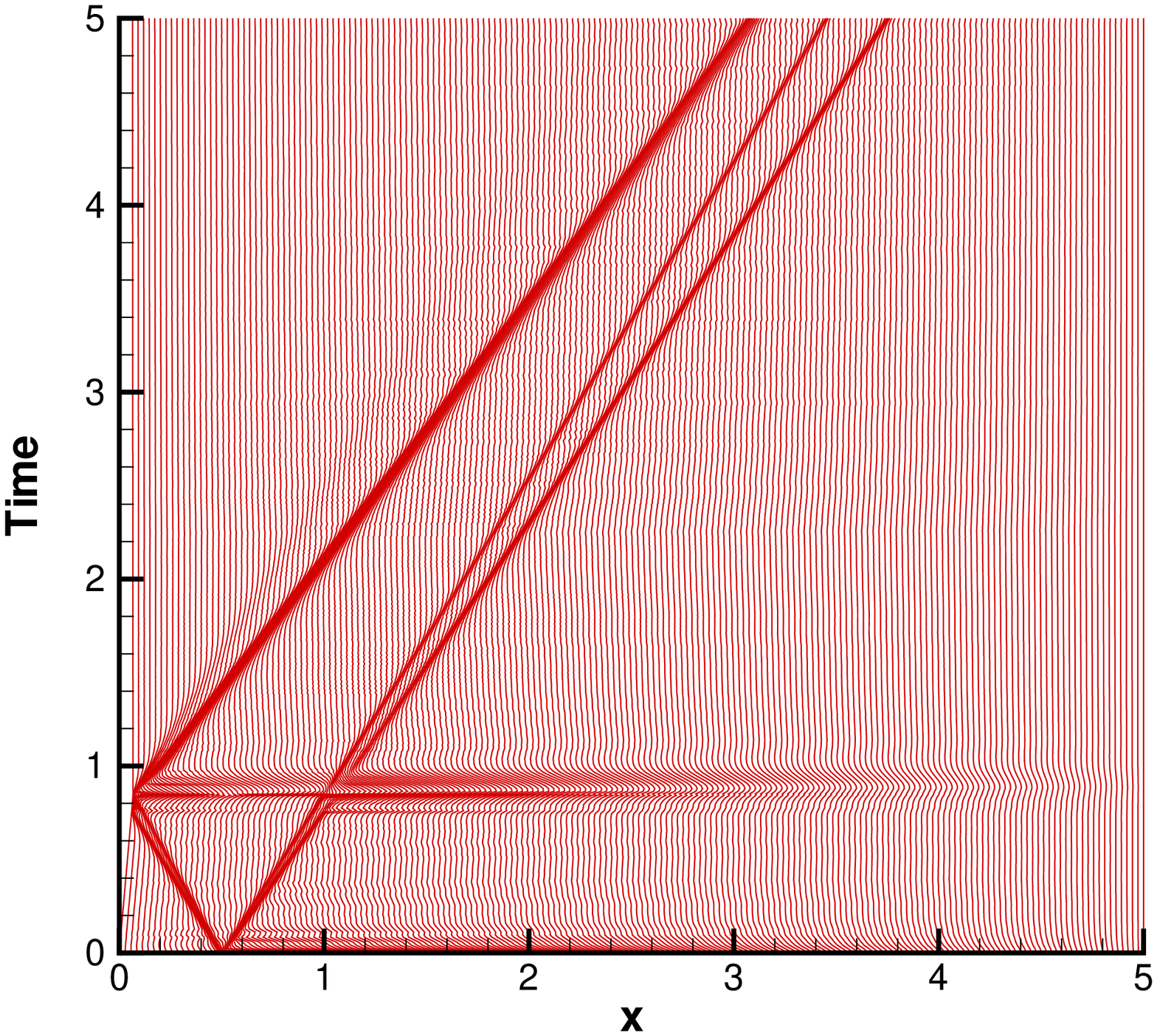}
\end{center}
\caption{Wilkins' problem. The mesh trajectories for CCL (left) and MMCC (right) schemes with 400 cells.} \label{fig2-xt}
\end{figure}

We also compare the numerical results of MMCC and CCL in Fig.~\ref{compar2-d}.
One can see that MMCC is more accurate than CCL for both shock and rarefaction waves.
Interestingly, the improvements are significant for the shock waves (between $3.4 < x < 4$) and the reflected
elastic rarefaction wave (between $3 < x < 3.2$) whereas that for the plastic rarefaction wave is only marginal.
This is reflected in the mesh concentration; see the mesh trajectories in Fig. \ref{fig2-xt}.
To explain this, we recall that the monitor function (\ref{mf2}) is defined in terms of the first derivatives
of $\rho$ and $s_{xx}$. The density changes abruptly in the shock wave regions while varying
only gradually in the rarefaction wave regions. Moreover, the flow in the regions of
both plastic shock and rarefaction waves is plastic, and from the von Mises yielding condition (cf. (\ref{e11})
and (\ref{e64})), $s_{xx}$ remains constant and its derivative is zero in these regions.
Consequently, the derivative of $s_{xx}$ is more significant in the regions of elastic (shock and rarefaction) waves
than in the regions of plastic waves. By combining the behavior of both $\rho$ and $s_{xx}$ and from
the definition of the monitor function, we know that the monitor function and therefore the mesh concentration
are large around the (elastic and plastic) shock waves and the elastic rarefaction waves.
On the other hand, the current definition of the monitor function leads to a relatively low mesh concentration
and only a marginal improvement in accuracy with MMCC over CCL around the plastic rarefaction wave.
Unfortunately, it is unclear to the authors how to concentrate more mesh points in the regions of plastic
rarefaction waves in an automatic manner. This will be an interesting topic for future research.

Overall, we have seen that MMCC is more accurate than CCL and is effective in concentrating the mesh
points around the shock and elastic rarefaction waves. No mesh crossing has been experienced
for this example too.

\section{Conclusions}
\label{SEC:conclusions}

In the previous sections a third-order moving mesh cell-centered scheme has been developed for
the numerical solution of one-dimensional elastic-plastic flows with the Mie-Gr\"{u}neisen equation
of state, the Wilkins constitutive model, and the von Mises yielding criterion.
The scheme moves and adapts the mesh by combining the Lagrangian method
with the MMPDE moving mesh method. The goal is to take advantages of both methods to better resolve
shock and other types of waves and be more robust in preventing the mesh from crossing and tangling.
Moreover, the MMCC scheme treats the relative velocity of the flow with respect to the mesh
as constant in time between time steps. An advantage of this treatment is that
free boundaries can be approximated with high-order accuracy. This is different from many existing moving mesh
methods where the mesh velocity is treated as constant, which gives a second-order
approximation of free boundaries.
Furthermore, a time dependent scaling is used in the computation of the monitor function to avoid
possible sudden movement of the mesh points due to the creation or diminishing of shock and rarefaction
waves or the steepening of those waves.
Finally, the two-rarefaction Riemann solver with elastic wave (TRRSE) is used to evaluate the
Godunov values of the density, pressure, velocity, and deviatoric stress at cell interfaces and
no remapping is used.

The third-order accuracy of the scheme has been verified for a smooth problem with Wilkins' constitutive model.
The MMCC scheme has also been used for the piston and Wilkins' problems. The
numerical tests have demonstrated the convergence of the scheme and its ability to concentrate
mesh points around shock and elastic rarefaction waves while the obtained numerical results are in good agreement
with those in literature. Comparison studies have shown that
MMCC is more accurate in resolving shock waves and rarefaction waves
than the CCL (Lagrangian mesh) scheme. It is also worth mentioning that no mesh crossing has been
experienced in the computation for all three examples.

Finally, it is pointed out that the focus of the current work is on one-dimensional problems.
Nevertheless, since both the cell-centered Lagrangian and MMPDE moving mesh methods are known
to work well in multi-dimensions, we believe that the MMCC scheme should be able to extend to
multi-dimensional elastic-plastic flows without major modifications. Investigations for such an extension
are underway.

\vspace{20pt}

\textbf{Acknowledgments.}
This work was partially supported by NSFC (China) (Grant No. 11172050, 11672047, and 11272064),
the National Basic Research Program (China) (Grant No. 2014CB745002),
and Science Challenge Project (Grant No. TZ2016002).
The authors are grateful to the anonymous referees for their valuable comments
in improving the quality of the paper.


\end{document}